\documentclass[11pt]{article}
\usepackage{fullpage, amsmath, amssymb, amsthm, epsfig, color, graphicx}
\newenvironment{@abssec}[1]{%
     \if@twocolumn
       \section*{#1}%
     \else
       \vspace{.05in}\footnotesize
       \parindent .2in
         {\bfseries #1. }\ignorespaces
     \fi}
     {\if@twocolumn\else\par\vspace{.1in}\fi}
\newenvironment{keywords}{\begin{@abssec}{Key words}}{\end{@abssec}}
\newenvironment{AMS}{\begin{@abssec}{AMS subject classification}}{\end{@abssec}}
\newtheorem{proposition}{Proposition}
\newtheorem{theorem}[proposition]{Theorem}
\newtheorem{lemma}[proposition]{Lemma}
\newtheorem{corollary}[proposition]{Corollary}
\newtheorem{result}[proposition]{Result}

\newtheorem{definition}[proposition]{Definition}
\newtheorem{claim}[proposition]{Claim}
\newtheorem{remark}[proposition]{Remark}
\newtheorem{conjecture}[proposition]{Conjecture}

\def\le{\leq}

\def\eqbd{\mathop{{:}{=}}}
\def\bdeq{\mathop{{=}{:}}}

\def\openC{{\rm C\kern-.48em\vrule width.06em height.6em depth-.02em
                 \kern.48em}}
\def\openQ{{{\rm Q\kern-.21cm\vrule width.6pt height 6.2ptdepth-.2pt \kern.21cm}}}
\def\openR{{{\rm I}\kern-.16em {\rm R}}}
\def\openZ{{{\rm Z}\kern-.28em{\rm Z}}}
\def\openT{{{\rm T}\kern-.42em {\rm T}}}
\def\openH{{{\rm I}\kern-.16em {\rm H}}}
\def\openK{{{\rm I}\kern-.16em {\rm K}}}
\def\openL{{{\rm I}\kern-.16em {\rm L}}}
\def\openM{{{\rm I}\kern-.16em {\rm M}}}
\def\openN{{{\rm I}\kern-.16em {\rm N}}}
\def\openP{{{\rm I}\kern-.16em {\rm P}}}

\def\eqbd{\mathop{{:}{=}}}

\def\ee{{\rm e}}
\def\dd{\, {\rm d}}

\let\N\openN

\let\R\openR
\def\Z{\mathbb Z}
\let\P\openP
\let\Q\openQ
\def\proof{\noindent {\bf Proof. \/}}
\def\proofof#1{\noindent {\bf Proof of #1. \/}}
\def\eop{\hfill
        {\ \vbox{\hrule\hbox{\vrule height1.3ex\hskip0.8ex\vrule}\hrule}}
        \vskip 0.3cm \par}

\def\spa{\mathop{\rm span}\nolimits}

\def\const{\mathop{\rm const}\nolimits}

\def\belowrightarrow#1{{{{}\over\ #1\ }\kern-1.1em\to}}
\def\l2{{L_2}}

\def\aint{\lceil \alpha \rceil{-}1}

\def\B{{\cal B}^+}

\begin{document}

\title{New coins from old, smoothly}
\author{Olga Holtz\thanks{Departments of Mathematics,
 University of California-Berkeley and Technische Universit\"at Berlin.
Research of O. Holtz was supported in part by a Center of Pure  and
Applied Mathematics grant at UC Berkeley and by the Sofja Kovalevskaja
Research Prize of the Humboldt Foundation, Germany.}
\and Fedor Nazarov\thanks{Department of Mathematics, University of
Wisconsin-Madison. Research of F. Nazarov was supported in part by NSF grants DMS-0501067 and DMS-0800243.}
\and Yuval Peres\thanks{
  Microsoft Research, Redmond. Research of Y. Peres was supported in part by NSF grant DMS-0605166.}}

\date{\small January 20, 2010}
\maketitle

\begin{keywords} Simulation,
approximation order, positive approximation,  Bernstein
operator, Lorentz operators, polynomial reproduction,  smoothness, H\"older class.
\end{keywords}

\begin{AMS} 41A10, 41A25, 65C50, 41A17, 68U20, 41A35, 41A27.
\end{AMS}

\begin{abstract}
Given a (known) function $f:[0,1] \to (0,1)$, we
consider the problem of simulating a coin with probability of heads $f(p)$ by tossing a coin
with unknown heads probability $p$, as well as a fair coin, $N$ times each, where $N$ may be random. The work of Keane and O'Brien (1994) implies that such a simulation scheme with the probability $\P_p(N<\infty)$
equal to $1$ exists iff $f$ is continuous. Nacu and Peres (2005) proved that $f$ is real analytic in 
an open set $S \subset (0,1)$ iff such a simulation scheme exists with the probability $\P_p(N>n)$ decaying 
exponentially in $n$ for every $p \in S$.
We prove that for $\alpha>0$ non-integer, $f$ is in the space $C^\alpha [0,1]$   if and only if a simulation scheme as above exists with
$\P_p(N>n) \le C (\Delta_n(p))^\alpha$, where  $\Delta_n(x)\eqbd \max \{\sqrt{x(1-x)/n},1/n \}$. 
The key to the proof is a new result in approximation theory:
 Let $\B_n$ be the cone of univariate polynomials 
with nonnegative Bernstein coefficients of degree $n$.
We show that a function  $f:[0,1] \to (0,1)$ is in $C^\alpha [0,1]$   if and only if
$f$ has a series representation $\sum_{n=1}^\infty F_n$ with $F_n \in \B_n$ and 
$\sum_{k>n} F_k(x) \le C(\Delta_n(x))^\alpha$ for all $ x \in [0,1]$ and $n \ge 1$. 
We also provide a counterexample to a theorem stated without proof by Lorentz (1963), 
who claimed that if some $\varphi_n \in \B_n$ satisfy $|f(x)-\varphi_n(x)|
\le C (\Delta_n(x))^\alpha$ for all $ x \in [0,1]$ and $n \ge 1$,
then $f \in C^\alpha [0,1]$.
\end{abstract}

\section{Introduction} Given a coin with unknown probability of heads $p\in [0,1]$, 
as well as a fair coin, we would like to simulate a coin with probability of heads 
$f(p)$ where $f:[0,1] \to (0,1)$ is a known function.
This means that we are allowed to toss the original $p$-coin and the fair coin $N$ times each, 
where $N$ is an almost surely finite stopping time (a notion to be clarified momentarily) 
 and declare heads or tails, depending on the 
outcome of these  $2N$ independent coin tosses. The probability of declaring a head must be exactly $f(p)$.

The measure corresponding to tosses of the $p$-coin is  the infinite product measure
 $\P_p$ on $\Omega=\{0,1\}^\N$ where in each coordinate the weights $(1-p,p)$ are used.
A measurable function $N:\Omega \to \N \cup \{\infty\}$ is a {\bf\em stopping time\/} 
if for every $k \in \N$, the indicator of $N=k$ is a function of the first $k$ coordinates in $\Omega$.
We say that $N$ is {\bf\em almost surely finite\/} if the probability $\P_p(N<\infty)$ is  $1$. More 
details on these notions can be found in any graduate textbook in Probability Theory, e.g.~\cite{Williams}.

This type of problem goes back to von Neumann's article~\cite{vonNeumann}
where he showed how to simulate a fair coin (i.e., $f(p)=1/2$) using only a biased $p$-coin where $p \in (0,1)$.
Moreover, the number of tosses $N$ needed satisfies $\P_p(N>n) \le \Bigl(1-2\epsilon(1-\epsilon)\Bigr)^{\lfloor n/2 \rfloor}$
if $p \in [\epsilon,1-\epsilon]$. In this paper we include a fair coin in the simulations since we want to consider $p$ near the endpoints $\{0,1\}$
where simulating a fair coin using a $p$ coin would be slow.

Since von Neumann's article, the simulation problem was subsequently solved for various
sother classes of functions -- see~\cite{KeaneOBrien, Peres, MosselPeres, NacuPeres}.
In particular, it was shown in~\cite{KeaneOBrien}  that an $f(p)$-coin can be simulated
using finitely many tosses of a $p$-coin for all $p$ in a closed interval $D\subseteq (0,1)$
if and only if $f$ is  continuous in $D$.
In~\cite{MosselPeres}, it was shown that  for $f:[0,1] \to (0,1)$, an $f(p)$-coin can be  simulated using finitely many tosses of  a $p$-coin via a
finite automaton for all $p\in (0,1)$, if an only if $f$ is a rational function over $\Q$.
(Simulation via a finite automaton is explained in detail in~\cite{MosselPeres}. An automaton 
is determined by a finite state space, a finite input alphabet, a transition rule from current 
state and input symbol to the next state, and a subset of final states. In our context, there 
are two final states, denoted $0$ and $1$, and we require that when the automation is given 
independent tosses of a $p$-coin as input, it will reach a final state with probability one, 
and output $1$ with probability $f(p)$.)

In \cite{NacuPeres}, it was shown that   if $D \subset (0,1)$ is closed and $f$ is real-analytic in an open neighborhood of $D$, then there is
a simulation of an $f(p)$-coin using $N$  tosses of  a $p$-coin where $N$ has uniform exponential tails for  $p \in D$, and conversely, if a simulation with exponential tails exists   for $p$ in an open set $S \subset (0,1)$,
then $f$ is real analytic in $S$. Moreover, the problem of simulation was recast in~\cite{NacuPeres} as
an approximation problem, and the question of characterizing simulation rates for non-analytic functions was posed.

\begin{definition}
Given a simulation algorithm, its {\bf\em simulation rate\/} is the probability
$\P_p \, (N>n)$ that  the number of required  inputs exceeds $n$.
(Each input is a toss of a $p$ coin and a toss of a fair coin).
If a simulation algorithm with  $ \P_p\, (N>n)= O( \psi_n(p))$  exists, we say
that the function $f$  {\bf\em can be simulated at the rate } $\psi_n(p)$.
\end{definition}

The goal of this paper is to show that the simulation rate is determined by the smoothness
of the simulated function $f$.  Our main result is that for positive $\alpha \notin \N$, a function $f:[0,1] \to (0,1) $ is in the space $C^\alpha$ (defined by a H\"older condition of order $\alpha-r$ on the derivative of order $r\eqbd \lfloor \alpha \rfloor$) if and only if $f$ can be simulated at the
rate $\Delta_n(p)^\alpha$, where $\Delta_n(x)\eqbd \max \{\sqrt{x(1-x)/n},1/n \}$, see Theorem \ref{thm-main} below.

\section{Preliminaries and statement of results}  \label{sec_prelim}

We first recall relevant definitions and results from the literature on this
problem and from approximation theory.   Recall that the univariate
{\bf\em Bernstein polynomials\/} of degree $n$ (see, e.g.,~\cite{Lorentz}) are defined
as \begin{equation}
 x\mapsto  p_{nk}(x)\eqbd {n\choose k} x^k (1-x)^{n-k}, \qquad k=0, \ldots, n. \label{B_basis}
\end{equation}
The Bernstein polynomials of degree $n$ form a basis for the space $\Pi_n$
of all polynomials of degree at most $n$. Thus, any polynomial $q$ of degree 
at most $n$ can be written as
$$ q(x)=\sum_{k=0}^n a_k p_{n,k}(x) ,   $$
 with the sequence $(a_0,\ldots,a_n)$
{\bf\em the degree $n$ Bernstein coefficients of\/} $q$. Whenever we write 
$$q\in {\cal B}_n,$$ 
this indicates that $q$ is already represented as a linear combination  of the Bernstein 
polynomials of degree $n$; this is admittedly an abuse of notation since the meaning 
of ``$q\in {\cal B}_n$'' differs from that of ``$q\in \Pi_n$''.
In addition, we write $$q\in \B_n$$ whenever the degree $n$ Bernstein coefficients 
$(a_0, \ldots, a_n)$ of $q$ are nonnegative. We will also need the following partial order 
on the space $\Pi_n$:
\begin{definition} \label{def_order}
Given  $q$, $r \in \Pi_n$, we write  $q\preceq_n r$, or $r\succeq_n q$, 
to denote that $r-q\in \B_n$. 
\end{definition}

Result~\ref{res_reduction} below was established in~\cite{NacuPeres} using
a simple probabilistic construction. This result reduces
the original simulation question to a problem in approximation theory, which we address
in this paper. In~\cite{NacuPeres} the goal was to obtain a simulation for $p$
in a closed subset of $(0,1)$; in this case a fair coin is not needed, as it can 
be produced from the $p$ coin using the von Neumann algorithm.
In the present paper we allow $p$ to range in the whole interval $[0,1]$,
so we use a fair coin in addition to the unknown $p$-coin.

\begin{result}[\cite{NacuPeres}] \label{res_reduction}
If there exists an algorithm that simulates a function $f$ on
a set $D\subset [0,1]$ using a random finite number $N$ of tosses of a
$p$-coin, then for all $n\geq 1$ there exist univariate polynomials
\begin{equation} \label{ghdef}
 g_n(x)\eqbd \sum_{k=0}^n {n\choose k}  a(n,k) x^k (1-x)^{n-k}, \qquad
h_n(x) \eqbd \sum_{k=0}^n  {n \choose k} b(n,k) x^k (1-x)^{n-k}
\end{equation}
with the following properties:
\begin{description}
\item[(i)] $0\leq a(n,k) \leq b(n,k) \leq 1$;
\item[(ii)] ${n \choose k} a(n,k)$ and ${n\choose k} b(n,k)$ are integers;
\item[(iii)] $\ g_n(p) \le f(p) \le h_n(p)$;
\item[(iv)] for all $m<n$ we have $g_m \preceq_n g_n$ and
$h_m \succeq_n h_n$;
\item[(v)] $h_n(p)-g_n(p)= \P_p \, (N>n)$.
\end{description}
Conversely, if there exist  polynomials $g_n$, $h_n$ as in  (\ref{ghdef}) satisfying
(i) -- (iv) with $\lim_n h_n(p)-g_n(p)=0$  for all $p\in D$, then there exists an algorithm that simulates an $f(p)$-coin for all $p \in D$   using
 $N$   tosses of the $p$-coin, where the random time $N$ satisfies $$ \P_p \, (N>n) = h_n(p)-g_n(p). $$

\end{result}

\noindent As noted in~\cite{NacuPeres}, given polynomials $g_n, h_n$ that satisfy all the requirements
except (ii), one can always round the values ${n\choose k} a(n,k)$ down and
the values ${n\choose k} b(n,k)$ up to an integer, and the resulting increase in the gap $h_n(p)-g_n(p)$ is exponentially small
in $n$ provided that  $p \in [\epsilon,1-\epsilon]$ for some $\epsilon>0$.
In the setting of the present paper, when the $p$-coin is tossed $n$ times we also
 toss a fair coin $n$ times; this means that
condition (ii) above is replaced by
\begin{itemize}
\item[{ \bf (ii')}] ${n \choose k} a(n,k)$ and $ {n\choose k} b(n,k)$ are integer multiples of $2^{-n}$,
\end{itemize}
since the probabilities of events that can be generated by tossing a fair coin $n$ times are precisely the integer multiples of $2^{-n}$.
Thus, given polynomials $g_n, h_n$ that satisfy  requirements (i), (iii) and (iv), we can round the  values ${n\choose k} a(n,k)$ down and
the values ${n\choose k} b(n,k)$ up to the nearest multiple of $2^{-n}$; this will only add at most $2^{1-n}$ to the gap $h_n-g_n$.

Therefore (up to an additional error term of $2^{1-n}$), the problem of determining the rate of simulation in our setting is equivalent to the problem of
determining the order of two-sided approximation to  $f$,  by polynomials $g_n$, 
$h_n \in {\cal B}_n$ that satisfy  requirements (i), (iii) and (iv).  We will
refer to requirements (iv) as the {\bf\em consistency requirements,\/} to the 
approximation scheme $(g_n)$ as a {\bf\em Bernstein-positive consistent approximation 
from below,\/} and to the approximation scheme $(h_n)$  as a {\bf\em Bernstein-positive
 consistent approximation from above.\/}

Observe that a Bernstein-positive consistent approximation to a function $f$
 from below is equivalent to  a certain nonnegative series representation of $f$.
Here is a precise statement.

\begin{lemma}\label{lem_series}
Let $D \subset [0,1]$ and let $(\psi_n)$ be a nonincreasing sequence of positive 
functions on $D$ that converges uniformly  to $0$.
A function $f$ is approximable on $D$ by a sequence of
Bernstein-nonnegative polynomials $(g_n)$ of degree $n$ satisfying the
consistency requirement (iv)
\begin{equation}   g_m  \preceq_n g_n \quad \hbox{\rm for all} \;\; n\geq m  \label{gn_iv}
\end{equation}
and the estimate
\begin{equation}
 0\leq f(x)-g_n(x) \leq \psi_n(x) \quad \hbox{\rm for all } \; x \in D \label{gn_rate}
\end{equation}
 if and only if $f$ can be represented
as a series
\begin{equation}  f(x)= \sum_{n=0}^\infty  F_n(x) \quad {\rm where} \quad
\sum_{n>N} F_n(x) \leq \psi_N(x) \quad \hbox{\rm for all } \; x\in D,
\label{series}
\end{equation}
where each $F_n$ is a polynomial in Bernstein form of degree $n$
with nonnegative coefficients.
\end{lemma}

\proof  Given an approximation scheme $(g_n)$ as above, set
$ F_n(x) \eqbd g_{n}(x)- g_{n-1}(x)$ where the second term $g_{n-1}(x)$
is rewritten in Bernstein form of degree $n$ and where $g_0(x)\eqbd 0$.
The consistency requirement~(\ref{gn_iv}) then guarantees that the
Bernstein coefficients of $F_n$ are nonnegative, and the
sum $\sum_{n>N} F_n(x)$ telescopes into $f(x)-g_N(x)$, which is
bounded pointwise by $\psi_N(x)$ according to~(\ref{gn_rate}).

Conversely, given a series representation~(\ref{series}), let
$g_n(x)\eqbd \sum_{k\leq n} F_k(x)$. Since the difference $F_n(x):= g_n(x)-g_{n-1}(x)$
is a Bernstein polynomial with nonnegative coefficients,
 the polynomials $(g_n)$ satisfy the consistency
requirement~(\ref{gn_iv}). Also, $f(x)-g_n(x)=\sum_{k>n} F_k(x) \leq \psi_n(x)$
 due to the rate condition~(\ref{gn_rate}). \eop

This approximation problem can be contrasted with the classical approximation
of a given function by (unrestricted) polynomials of degree at most $n$ on the interval $[0,1]$.
In that case, the approximation order coincides with the smoothness of $f$.
To state this classical result precisely, we first recall how smoothness is
measured.

\begin{definition} Let $\alpha>0$ with $\alpha \notin \N$. A function $f$ is said to be in the {\bf\em smoothness class\/}
$C^\alpha[0,1]$  if $f$ is  $r\eqbd \lfloor \alpha \rfloor$  times differentiable
and the following condition holds:

The {\bf\em modulus of continuity\/} of the
$r$th derivative $f^{(r)}$
$$ \omega(f^{(r)};h) \eqbd \sup_{x,y \in [0,1],\; |x-y|<h} |f^{(r)}(x)-f^{(r)}(y)| $$
is of order  $O(h^{\alpha-r})$. In that case, we will use the
notation $$ \|f \|_{C^\alpha} \eqbd \sup_{h>0}  {\omega(f^{(r)};h) \over h^{\alpha - r}  }.$$
(Note this is a seminorm rather than a norm, as it vanishes on  polynomials of degree at most $r$.)

\end{definition}


The order of approximation of a given function by polynomials is then determined as follows.

\begin{result}[{see, e.g.,~\cite[Chapter 8,~Theorem~6.3]{DeVoreL}}] \label{res_best}
Let $\alpha>0$ be a non-integer. There exists a sequence of polynomials $\{p_n\}$, where  
the degree of $p_n$ is at most $n$, satisfying
$$ | p_n(x) -  f(x)| = O((\Delta_n(x))^{2\alpha}) \qquad \hbox{\rm for all}
\;\; x\in [0,1] \, ,$$ if and only if $f\in C^\alpha[0,1]$.
Here the quantity $\Delta_n(x)$ is defined by
$$ \Delta_n\eqbd \Delta_n(x)\eqbd \max \left\{\sqrt{ x(1-x)\over n} \, , \,
 {1\over n} \right\}.$$
\end{result}

In other words, the rate of approximation of $f\in C^\alpha[0,1]$ is $O(n^{-\alpha})$
away from the boundary of the interval $[0,1]$ and is $O(n^{-2\alpha})$ close to the endpoints $0$ and $1$.
The characterization of the rate of polynomial approximation for integer values $\alpha$
involves the   generalized Zygmund class, which we will recall in
Section~\ref{sec_moments}. In the main part of this paper, we work under the assumption
$\alpha \notin \N$.

Result~\ref{res_best} shows that   a function $f\in C^\alpha[0,1] \setminus C^{\alpha+\epsilon}[0,1]$
 cannot be simulated at the rate $O(\Delta_n^{2(\alpha+\epsilon)})$. However, since our approximants must satisfy special
restrictions imposed by  Result~\ref{res_reduction}, we should not expect to achieve
the approximation order provided by unrestricted polynomials of degree $n$.


In view of requirement (i), it is natural to consider first the approximation order
achieved by polynomials with nonnegative Bernstein coefficients.
G.~G.~Lorentz proposed a solution to this problem in~\cite{Lorentz_paper},
where he argued that the approximation order under this constraint is half the approximation
order provided by unconstrained polynomials, i.e., half the smoothness of the function $f$.
Theorem~1 of \cite{Lorentz_paper} establishes that a $C^\alpha$-function $f$ can
be approximated at the rate $O((\Delta_n)^\alpha)$ by Bernstein-nonnegative polynomials
of degree $n$.

\begin{result}[{\cite[Theorem~1]{Lorentz_paper}}] \label{result-L}
Let $\alpha>0$. A positive function $f \in C^\alpha[0,1]$ can be approximated
by polynomials $q_n$ of degree at most $n$ with nonnegative Bernstein coefficients  at the rate
\begin{equation}  \label{claim-L}
 | q_n(x) -  f(x)| = O((\Delta_n(x))^\alpha) \qquad
\hbox{\rm for all} \;\; x\in [0,1].
\end{equation}
\end{result}
Lorentz~\cite{Lorentz_paper} also stated (without proof) a converse to this result; unfortunately, that converse is incorrect.
 We return to this point at the end of the section.

We will  use a variant of Lorentz' approach to establish our main result, that
with   the extra requirements (i), (iii), (iv) in place, we can still achieve the
same approximation order as in~(\ref{claim-L}).
\begin{theorem}  \label{thm-main}
Let $f: [0,1] \to (0,1)$ and  let $\alpha>0$ with $\alpha \notin \N$. If $f \in C^\alpha [0,1]$, then $f$ can be simulated at the rate $(\Delta_n(x))^\alpha$ on $[0,1]$. Precisely,
there exist polynomials $g_n$ and $h_n$ satisfying conditions (i), (ii'), (iii) and (iv) of
Result~\ref{res_reduction} and
 $h_n(x)-g_n(x) = O( (\Delta_n(x))^\alpha )$  uniformly in $[0,1]$.
Conversely, if $f$ can be simulated  at the rate $(\Delta_n(x))^\alpha$ on the interval
$[0,1]$,  then $f\in C^\alpha [0,1]$.
\end{theorem}

\medskip

We begin by proving a reduction lemma that shows that it is enough to find
consistent approximants $g_n$, $h_n$ for each $b$-adic degree $n \in b^{\N}=\{b^\ell\}_{\ell \ge 1}$ where
$b$ is a fixed integer greater than $1$. Using
these approximants, one can then interpolate between $b$-adic levels to build
up a consistent approximation scheme providing the same approximation order as
the $b$-adic polynomials $g_n$, $h_n$, $n=b^\ell$. This $b$-adic idea {\em per se\/}
is quite well known  and, in particular, is used in~\cite{NacuPeres} with $b=2$.

\begin{lemma} \label{lem_reduce} Let $b$ be a fixed integer greater than $1$.
Given a function $f$ on $[0,1]$, suppose there exist two sequences of
polynomials  $\left(g_n\right)_{n \in b^\N}$, $\left(h_n\right)_{n \in b^\N}$   satisfying  conditions~(i), (iii), (iv)  of Result~\ref{res_reduction}
with $\psi_n$ of order $O((\Delta_n)^\alpha)$, so that
\begin{equation}
 h_n(x)-g_n(x) = O( (\Delta_n(x))^\alpha ) \qquad \hbox{\rm uniformly in }   [0,1]
 \label{app_order}
\end{equation}
for $n \in b^\N$.
Then these sequences can be augmented to full sequences $\left(g_n\right)_{n\in \Z_+}$,
$\left(h_n\right)_{n\in\Z_+}$  satisfying conditions (i), (iii), (iv) from
Result~\ref{res_reduction} and  condition~(\ref{app_order}) for all $n\in Z_+$.
In particular, under these assumptions there exists an algorithm that simulates an 
$f(p)$-coin at the rate
$O((\Delta_n(p))^\alpha)$ on $[0,1]$.
\end{lemma}

\proof Given the polynomials $g_n$ and $h_n$ for $b$-adic values of $n$,
we will fill in the gaps in the two sequences in the obvious way:
given $n$, let $n'\eqbd b^{\lfloor \log_b n \rfloor}$, and set
$$ g_n(x) \eqbd (x+(1-x))^{n-n'} g_{n'}(x,y) , \qquad
h_n(x) \eqbd (x+(1-x))^{n-n'} h_{n'}(x,y) $$
by expanding the right-hand sides into Bernstein polynomials of degree $n$.
The Bernstein coefficients of the resulting polynomials $g_n$, $h_n$
are therefore some convex combinations of the coefficients of $g_{n'}$,
$h_{n'}$. It follows that condition (i) holds for the full sequences
$(g_n)$, $(h_n)$. It is clear from the construction that (iii) and (iv)
hold as well, the latter condition being an equality
except when jumping from one $b$-adic level to the next, when it is
satisfied by our assumption. Recall that condition (ii') can always be
satisfied by introducing an exponentially small correction, so there is no need to verify it explicitly.
To check that~(\ref{app_order}) holds for the full sequences $(g_n)$, $(h_n)$,
note that, by construction,
$$ h_n(x)-g_n(x)=h_{n'}(x)-g_{n'}(x)=O((\Delta_{n'}(x))^\alpha).$$
But since $n' \leq n < b n'$, we see that $O((\Delta_{n'}(x))^\alpha)=O((\Delta_n(x))^\alpha)$.
This completes the proof.    \eop

As noted already, Lorentz~\cite{Lorentz_paper} stated a converse to Result \ref{result-L} above, which (in a special case) can be written as follows.

\begin{claim}[{\cite[Theorems~5 and 6]{Lorentz_paper}}]   \label{res_Lorentz}
Let $\alpha>0$. If a function $f$ can be approximated by polynomials $q_n$
of degree at most $n$ with nonnegative Bernstein coefficients  at the rate
(\ref{claim-L}), then $f\in C^\alpha[0,1]$.
\end{claim}
The  argument proposed in~\cite{Lorentz_paper} for these theorems
skips technical details and refers  to the work of Timan \cite{Timan_book}.
Specifically, we quote Theorem~6 from~\cite{Lorentz_paper} and the subsequent
discussion:

\begin{quote}
`` {\bf Theorem 6.} {\em For each $r=1,2,\ldots$ there is a constant $C_r$ with the following
property. Let $\omega(h)$ be a modulus of continuity, and put
$$ \tilde{\omega}(h)=h \int_h^1 {\omega(u)\over u^2} \dd u +
\int_0^h {\omega(u) \over u } \dd u .   $$
If $f(x)$ is a continuous function on $[0,1]$ and if there exists a sequence
$P_n(x)$ of polynomials with positive coefficients of degree $n$ such that
$$  |f(x)- P_n(x)| \leq (\Delta_n)^r \omega (\Delta_n), \quad  0\leq x \leq 1, \quad
n=0,1,\ldots ,  $$
then $f$ has on $[0,1]$ the continuous derivatives $f'$, $f''$, $\ldots$, $f^{(r)}$ and
$$  \omega(f^{(r)}; h) \leq C_r \tilde{\omega}(h). $$ }
  We omit the proofs. The method of deriving theorems of this kind from inequalities of the
Markov-Bernstein type is due essentially to S.~Bernstein, and is well known. For the
variation of it which fits the present situation especially well, compare [6, p.\ 357
and p.\ 360]. It should be noted that [6, p.\ 357] contains an essential mistake: the
derivative $P'_{2^{m+1}}$ on p.\ 359 should have been estimated at a point different
from $x$. However, the proof can be corrected. ''
\end{quote}

\noindent
In this quote, [6] refers to the original Russian  edition of Timan's work
\cite{Timan_book}.
Trying to reconstruct Lorentz' complete  argument for his Theorem 6, we came to the realization
that his argument requires an extra assumption, in fact precisely the
assumption of Bernstein-nonnegative consistent approximation, or equivalently,
the nonnegative series representation~(\ref{series})  that is central to this paper.
In the next section we show that, indeed, such a series representation of $f$
with tails decaying at the rate $(\Delta_n)^{\alpha}$ implies the
$C^\alpha$ smoothness of the represented function $f$.  Thus our results
here also provide a correction to the statement  of Lorentz.
In  Section~\ref{sec_fedya}, we construct a counterexample to
Theorem~6 from~\cite{Lorentz_paper}.

Our final point in this section concerns notation. In the rest of the paper, we will
establish a number of estimates on various functions. The constants
in such estimates will be usually simply denoted by $\const$ or, say, $\const_j$, the
latter indicating that the constant may depend on $j$. A few constants that are crucial
to our main argument will be labeled by the number of the theorem
or lemma where they occur.

\section{Consistent approximation implies smoothness in Theorem \ref{thm-main}} \label{sec-necessity}

Lorentz proved the following analogues of Bernstein's and Markov's inequalities
 (both original inequalities can be
found, e.g., in~\cite{DeVoreL}). This result of Lorentz is formulated
for a certain class of functions $\Omega$; we will use it only for the
power functions $t \mapsto t^j$.

\begin{result}[{\cite[Theorem 3]{Lorentz_paper}}]   \label{res_Lorentz_neg}
For each $r=1,2,\ldots$ and each $H>0$, there is a constant $K_r = K_r(H)$
with the following property. If $\Omega(h)$ is an increasing positive function
defined for all $h\geq 0$ such that
$$ \Omega(2h) \leq H \Omega(h), \quad h\geq 0,$$
then for each Bernstein-positive polynomial $P_n$ of
degree $n$, the inequality
$$ P_n(x) \leq \Omega(\Delta_n(x)), \quad 0\leq x \leq 1 $$
 implies
\begin{equation}
 |P^{(r)}_n(x)| \leq K_r (\Delta_n(x))^{-r}\, \Omega(\Delta_n(x)), \quad 0\leq x \leq 1.
\label{Markov-Bernstein}
\end{equation}
\end{result}

We need the following observation.

\begin{lemma}\label{lem-1} For any $x$ and $\xi$ in $[0,1]$,
\begin{align*}
\max \left\{ \frac{\Delta_n(\xi)}{\Delta_n(x)}, \frac{\Delta_n(x)}{\Delta_n(\xi)} \right\}
\leq 2\left(1+ \frac{|x- \xi|}{\Delta_n(x)}\right)~.
\end{align*}
\end{lemma}

\proof We start by proving one of the two inequalities, viz.,
 \begin{equation}
   \frac{\Delta_n(x)}{\Delta_n(\xi )} \le 2\Bigl(1+\frac{|x-\xi |}{\Delta_n(x)}\Bigr)\,.
\label{Delta-bound}
\end{equation}
By the symmetry $\Delta_n(x)=\Delta_n(1-x)$, we may assume that $x,\xi  \in [0,1/2]$. We also assume that $\xi <x$ and $\Delta_n(x)>1/n$, since otherwise the inequality is obvious. If $\xi  \ge x/2$ then the left-hand side of (\ref{Delta-bound}) is at most $2$, so we may assume that $\xi <x/2$. In this case we have
$$(\Delta_n(x))^2 \le x/n \le 2|x-\xi |/n \le 2|x-\xi |\, \Delta_n(\xi ) ,
$$ which implies (\ref{Delta-bound}).

The proof of the other inequality (which bounds  $\Delta_n(\xi)/\Delta_n(x)$ by the right-hand side of
(\ref{Delta-bound})) is very similar. We may again assume, by symmetry,
that  $x$, $\xi \in [0,1/2]$. We also assume that $\xi>x$ and $\Delta_n(\xi)>1/n$
since otherwise the inequality is obvious. Thus $\Delta_n(\xi) \leq \sqrt{\xi/n}  < \xi$.
If $\xi \leq 2x$, then the left-hand side is at most $2$, while the right-hand
side is greater than $2$. Thus, the only remaining case is $\xi>2x$. Then
$$ \Delta_n (\xi) \leq \xi \leq 2 |\xi - x|. $$
Dividing by $\Delta_n(x)$, we obtain the desired bound.
  \eop

To prove the necessity of $C^\alpha$-smoothness, we will follow the approach suggested by
G.~Lorentz in~\cite{Lorentz_paper}, which goes back to Timan~\cite{Timan_book}
and ultimately to S.~Bernstein.  \bigskip

\proofof{necessity in Theorem~\ref{thm-main}}  
Suppose that  $f$ can be simulated at the rate
$(\Delta_n)^{\alpha}$ on the interval $[0,1]$.
 Using the sequence $(g_n)$ that approximates
$f$ from below and satisfies the consistency requirement
$g_n \preceq_{2n} g_{2n}$,
we set $G_n\eqbd g_{2^{n+1}}-g_{2^n}$ and obtain the following nonnegative
series  representation for $f$:
\begin{equation}
 f(x) = \sum_{n=0}^\infty  G_n(x), \qquad
G_n \in \B_{2^n}. \label{2-series}
\end{equation}
By the assumption on the rate of approximation, the polynomials $G_n$
satisfy the bound $$  |G_n(x)| \leq \const (\Delta_{2^n}(x))^\alpha
\qquad
\hbox{\rm for all }\; x\in [0,1]. $$
Now, the inequality~(\ref{Markov-Bernstein}) implies
\begin{equation}
  |G^{(j)}_n(x)| \leq \const (\Delta_{2^n}(x))^{\alpha-j} \qquad \hbox{\rm for all }\;
x\in [0,1], \; j\in \N.  \label{f-Bernstein}
\end{equation}
This already ensures that we can differentiate (\ref{2-series}) term by term $r$ times, and that
\begin{equation}
 f^{(r)}(x) = \sum_{n=0}^\infty  G_n^{(r)}(x), \label{2r-series}
\end{equation}
is continuous in $[0,1]$.
Our goal is to prove that $f\in C^\alpha[0,1]$, i.e., that the inequality
\begin{equation} \label{hold}
  \left|  f^{(r)}(x) - f^{(r)}(y)\right| \leq \const |x-y|^{\alpha-r}
\end{equation}
holds for $x,y \in [0,1]$.  Without loss of generality $x(1-x)\ge y(1-y)$, whence
$\Delta_n(x) \ge \Delta_n(y)$ for all $n$.
For any $n$, there is some $\xi_n$ between $x$ and $y$  such that
\begin{equation} \label{meanval}   |G^{(r)}_n(x) - G^{(r)}_n(y)|=|x-y| \, |G_n^{(r+1)}(\xi_n)| \le \const |x-y|\,  (\Delta_{2^n}(\xi_n))^{\alpha-r-1} \,,
\end{equation}
using the bound~(\ref{f-Bernstein}) with $j=r+1$.
Choose $N$ so that
\begin{equation} \label{Nchoice}
\Delta_{2^{N+1}}(x) < |x-y| \leq \Delta_{2^{N}}(x) \,.
\end{equation}
For $n \le N$ we have $|x-y| \leq \Delta_{2^{n}}(x)$, so  Lemma \ref{lem-1} implies that
$\Delta_{2^n}(x) \le 4 \Delta_{2^n}(\xi_n)$. Thus for $n \le N$, (\ref{meanval}) gives
\begin{equation} \label{meanval2}   |G^{(r)}_n(x) - G^{(r)}_n(y)|  \le \const |x-y|\,  (\Delta_{2^n}(x))^{\alpha-r-1} \,.
\end{equation}
We now write $f^{(r)}(x)-f^{(r)}(y)$ by splitting the sum~(\ref{2r-series}) into two parts:
\begin{equation} \label{decomp} f^{(r)}(x) - f^{(r)}(y) = \sum_{n=0}^{N} ( G^{(r)}_n(x) - G^{(r)}_n(y))
 +   \sum_{n=N+1}^\infty  (G^{(r)}_n(x)-G^{(r)}_n(y)).
\end{equation}
Estimate the first sum using~(\ref{meanval2})
and the second using~(\ref{f-Bernstein}) with $j=r$. This yields
\begin{eqnarray}
 \left|  f^{(r)}(x) - f^{(r)}(y)\right|& \leq  & \const \Bigl( |x-y| \sum_{n=0}^{N} (\Delta_{2^n}(x))^{\alpha-r-1} + 2
\sum_{n=N+1}^\infty (\Delta_{2^n}(x))^{\alpha-r} \Bigr) \label{useybound} \\
& \leq & \const |x-y| (\Delta_{2^N}(x))^{\alpha-r-1} + \const (\Delta_{2^{N+1}}(x))^{\alpha-r} \,,
\label{last_estim}
\end{eqnarray}
where we used the inequality $\Delta_{k}(x) \ge \sqrt{2} \, \Delta_{2k}(x)$ to compare the two series in (\ref{useybound}) to geometric series.
In view of (\ref{Nchoice}), the bound~(\ref{last_estim}) yields (\ref{hold}).
 \eop

In the preceding proof, the strict inequality $\alpha<r+1$ is used only at one point:
to show that the sum of terms with $n \le N$ in (\ref{decomp}) is comparable to the last term.
(If $\alpha=r+1$ then all these terms are of the same magnitude and we lose a factor of
 $N \approx \log \frac{1}{|x-y|}$ in the estimate.) Nevertheless, for the case
 $\alpha=r+1$, the same method will allow us to show that $f^{(r)}$ is in the Zygmund class.

\begin{theorem}  \label{thm-int}
Let $r$ be a nonnegative integer. Suppose that $f:[0,1] \to (0,1)$  can be simulated at the rate 
$(\Delta_n(x))^{r+1}$ on $[0,1]$. Precisely, suppose that
there exist polynomials $g_n$ and $h_n$ satisfying conditions (i), (ii'), (iii) and (iv) of
Result~\ref{res_reduction} and  $h_n(x)-g_n(x) = O( \Delta^{r+1}_n(x) )$   uniformly in $[0,1]$. 
Then  $f^{(r)}$ is in the {\bf \em Zygmund class,\/} that is
\begin{equation} \label{rzyg}
 |f^{(r)}(x+\delta)-2f^{(r)}(x)+f^{(r)}(x-\delta)|=O(\delta) \,,
 \end{equation}
 uniformly for all $x, \delta$ such that $0 \le x-\delta< x+\delta \le 1$.
\end{theorem}
In fact, as in the preceding proof, only the approximation from below by $g_n$ is used.
\medskip

\proof  The hypothesis implies that $f$  has a
series  representation  as in (\ref{2-series})
where the polynomials $G_n \in \spa_+ \{p_{2^n, k}: k=0, \ldots, 2^n  \}$
satisfy  $  |G_n(x)| \leq \const (\Delta_{2^n}(x))^{r+1}$ for all $x\in [0,1]$.
The inequality~(\ref{Markov-Bernstein}) implies
\begin{equation}
  |G^{(j)}_n(x)| \leq \const (\Delta_{2^n}(x))^{r+1-j} \qquad \hbox{\rm for all }\;
x\in [0,1], \; j\in \N \, , \label{f2-Bernstein}
\end{equation}
so  (\ref{2r-series}) holds and $f^{(r)}$
is continuous in $[0,1]$. Fix $\delta \in (0,1/2)$, and choose $N$ minimal so that 
$\delta \leq \Delta_{2^{N}}(x) \,.$
Write $f=S_1+S_2$ where
$$ S_1= \sum_{n=0}^N  G_n(x) \quad \mbox{ \rm and } S_2 = \sum_{n=N+1}^\infty  G_n(x) \,. $$
The preceding proof works to show   that $S_2(x+\delta)-S_2(x)=O(\delta)$, and this implies that the estimate
(\ref{rzyg}) holds with $S_2^{(r)}$ in place of $f^{(r)}$.  It remains to handle $S_1^{(r)}$.

For any $n$, there is some $\eta_n \in [x-\delta,x+\delta]$   such that
\begin{equation} \label{meanval3}   |G^{(r)}_n(x+\delta) - 2G^{(r)}_n(x)+ G^{(r)}_n(x-\delta) |=\delta^2 \, |G_n^{(r+2)}(\eta_n)| \le \const \delta^2\,  (\Delta_{2^n}(\eta_n))^{-1} \,,
\end{equation}
using the bound~(\ref{f2-Bernstein}) with $j=r+2$.

For $n \le N$ we have $\delta \leq \Delta_{2^{n}}(x)$, so  Lemma \ref{lem-1} implies that
$\Delta_{2^n}(x) \le 4 \Delta_{2^n}(\eta_n)$. Thus for $n \le N$, (\ref{meanval3}) gives
$$ |G^{(r)}_n(x+\delta) - 2G^{(r)}_n(x)+ G^{(r)}_n(x-\delta) | 
 \le \const \delta^2 (\Delta_{2^n}(x))^{-1} \,. $$
This yields
\begin{eqnarray}
 \left|  S_1^{(r)}(x+\delta) - 2S_1^{(r)}(x) + S_1^{(r)}(x+\delta)\right|& \leq  & \const
 \Bigl( \delta^2 \sum_{n=0}^{N} (\Delta_{2^n}(x))^{-1}  \Bigr) \nonumber \\
& \leq & \const \delta^2   (\Delta_{2^N}(x))^{-1} \, \le \const \delta \,.
\label{last_estim2}
\end{eqnarray}
 The previous estimate for $S_2$, together with the bound~(\ref{last_estim2}), yields (\ref{rzyg}).
 \eop

\section{Lorentz operators and simultaneous approximation} \label{sec-operators}
In the following three sections, we restrict attention to $\alpha \notin \N$.
We now recall the main ingredients of the valid proof of Result \ref{result-L} (Theorem~1
from \cite{Lorentz_paper}). That proof is based on the Taylor expansion
\begin{equation}
f(x) = f\left({k\over n}\right) -\sum_{j=1}^r {1\over j!} \left( {k\over n} -x \right)^j f^{(j)}(x)
+{1\over r!} \left({k\over n} -x \right)^r [f^{(r)}(x)-f^{(r)}(\xi_k)],
\label{Taylor}
\end{equation}
where $\xi_k \eqbd \xi_k(x)$ is a point between $x$ and $k/n$ and $f$ is assumed
to be  $r$ times differentiable. This formula is used in~\cite{Lorentz_paper} to
derive an asymptotic expansion of the {\bf\em Bernstein operator\/}
$$ (B_n f)(x)\eqbd \sum_{k=0}^n  f\left({k\over n}\right) p_{nk}(x),  $$
where the polynomials $p_{nk}$ are defined in~(\ref{B_basis}).
Multiplying the Taylor expansion~(\ref{Taylor}) by $p_{nk}(x)$ and summing over $k$,
we obtain
\begin{eqnarray}
f(x) & = & (B_n f)(x) -\sum_{j=1}^r {1\over j! n^j} T_{nj}(x)
 f^{(j)}(x)+ (R_r f) (x), \quad {\rm where} \label{f_recur}  \\
T_{nj}(x) & \eqbd & \sum_{k=0}^n (k-nx)^j p_{nk}(x),  \label{Ts} \\
(R_r f) (x) & \eqbd & {1\over r!} \sum_{k=0}^{n} \left({k\over n} -x \right)^r
 [f^{(r)}(x)-f^{(r)}(\xi_k)] p_{nk}(x). \nonumber
\end{eqnarray}
This leads Lorentz to the natural definition of the operators $Q_{n,r}$, using
 the recurrence
\begin{equation}
(Q_{n,r}f)(x)  \eqbd  (B_n f)(x) -\sum_{j=1}^r {1\over j! n^j} T_{nj}(x)
(Q_{n,r-j} f^{(j)})(x),   \label{Qs}
\end{equation}
where each $f^{(j)}$ in~(\ref{f_recur}) is replaced by its approximation
$Q_{n,r-j} f^{(j)}$.

Note that the sum in~(\ref{Qs}) in
fact starts at $j{=}2$ rather than at $j{=}1$, since the polynomial $T_{n1}$ is identically
zero. Also note that the expressions~(\ref{Qs}) must be written in the Bernstein basis
of degree $n{+}r$, so that, e.g., the leading term $B_n f$
must be multiplied by the binomial expansion of $(x+(1-x))^r$ to appear in its
Bernstein form of degree $n{+}r$. We will refer to the operators $Q_{n,r}$
mapping a function to a polynomial in Bernstein form of degree $n{+}r$ as
the {\bf\em Lorentz operators.\/}

An important property of the Lorentz operators that can be inferred
directly from their recursive definition is their {\bf\em polynomial reproduction.\/}
Precisely,  the Lorentz operator  $Q_{n,r}$ reproduces polynomials
of degree at most $r$.

\begin{lemma}  \label{lem_rep}
Let $f$ be a polynomial of degree at most $r$. Then
$Q_{n,r}f=f$.
\end{lemma}

\proof The proof is by induction on $r$.
The result holds for $r=0$ and $1$ since $Q_{n,0}=Q_{n,1}$ is simply
the Bernstein operator, which reproduces linear functions.
For higher values of $r$, the proof is as follows.
The Taylor polynomial of $f$ of degree $r$ coincides with $f$, so
$$ f(x) = f\left({k\over n}\right) -\sum_{j=1}^r {1\over j!}
\left( {k\over n} - x \right)^j f^{(j)}(x),  $$
so by multiplying by $p_{nk}(x)$ and summing over $k$, we obtain
$$ f(x) = (B_n f)(x) -\sum_{j=1}^r {1\over j! n^j} T_{nj}(x) f^{(j)}(x).  $$
By our inductive assumption, $f^{(j)}=Q_{n,r-j}f^{(j)}$. Substituting
this into~(\ref{Qs}), we get $Q_{n,r}f=f$.
\eop

 As noted in \cite{Lorentz_paper}, the Lorentz operators  can be rewritten as follows
$$ (Q_{n,r} f)(x) \bdeq \sum_{k=0}^n \left( f\left({k\over n} \right) +\sum_{j=2}^r
f^{(j)} \left( {k\over n}  \right) {1\over n^j}\, \tau_{rj} (x,n)  \right) p_{nk}(x), $$
 or, more simply, as
\begin{equation}
(Q_{n,r}f)(x)  = \sum_{k=0}^n \left( \sum_{j=0}^r f^{(j)}\left(  {k\over n} \right)
{1\over n^j} \tau_{rj}(x,n) \right)  p_{nk}(x)  \label{simpleQs},
\end{equation}
with the understanding that $\tau_{r0}(x,n)=1$ and $\tau_{r1}(x,n)=0$.
Plugging~(\ref{simpleQs}) into~(\ref{Qs}), we obtain
$$ (Q_{n,r}f)(x)=\sum_{k=0}^n \left(  f\left({k\over n} \right) -\sum_{j=2}^r
f^{(j)} \left( {k\over n} \right) {1\over n^j} \sum_{l=2}^j {1\over l!}
T_{nl}(x,n) \tau_{r-l,j-l}(x,n)   \right) p_{nk}(x). $$
By term-by-term comparison, this yields
\begin{equation}
 \tau_{rj}(x,n)=-\sum_{l=2}^j {1\over l!} T_{nl}(x) \tau_{r-l,j-l}(x,n) \qquad
{\rm for} \;\; j\geq 2. \label{recur_taus}
\end{equation}
The recurrence~(\ref{recur_taus}) can be used to show that the polynomials
$\tau_{rj}(x,n)$ are independent of $f$, are of degree $j$ in $x$ and of
degree $\lfloor j/2 \rfloor$ in $n$, as  noted by Lorentz \cite{Lorentz_paper}.
The recurrence~(\ref{recur_taus}) also shows that, as functions,
the $\tau_{rj}$  do not depend on the parameter $r$. However, in the
expression for $Q_{n,r}$, each of the  $\tau_{rj}$s is written in its Bernstein
form of degree $r$ to bring the entire expression $Q_{n,r} f$ into its Bernstein form
of degree $n+r$. Since we are mainly concerned with pointwise estimates on the
$\tau_{rj}$s, we will use the simpler notation   $\tau_j \eqbd \tau_{rj}$.
We will begin with the following estimates on the polynomials $\tau_{j}$:

\begin{lemma} \label{lem_taus} The polynomials $\tau_{j}$ are bounded  by
\begin{equation}
| \tau_{j}(x,n)  |  \leq \const_j \,  n^j (\Delta_n(x))^j \qquad \hbox{\rm for all} \;\;
x\in [0,1],\label{bound_tau}
\end{equation}
where $\const_j$ depends only on $j$.
\end{lemma}

\proof We use induction on $j$. For $j=0,1$ (\ref{bound_tau}) is clear.
 By~\cite[(17) on p.~242]{Lorentz_paper},
$$ | T_{n\ell}(x) |  \leq  \const_\ell\, n^l (\Delta_n(x))^\ell \,.$$
Applying~(\ref{recur_taus}) and the induction hypothesis
$$\forall \ell>0, \quad  | \tau_{j-\ell}(x,n)| \leq \const_{j-\ell} \, n^{j-\ell} 
(\Delta_n(x))^{j-\ell}$$
  gives~(\ref{bound_tau}), as required.  \eop

\begin{corollary}  \label{cor_taus}
Fix an integer $r \ge 0$. For any $j \le r$, write
$$ \tau_{j}(x,n)\bdeq \sum_{i=0}^j a_i(n,j)x^i(1-x)^{j-i}\,. $$
Then for all $i \in [0,j]$, we have $|a_i(n,j)|=|a_{j-i}(n,j)|$ and $|a_i(n,j)| \le  C_j^\sharp n^{i }$
for some constants $C_j^\sharp$.
\end{corollary}

\proof The polynomials $T_{nj}$ satisfy $T_{nj}(1-x)=(-1)^jT_{nj}(x)$, as is easily seen using  the substitution $\tilde{k}=n-k$ in their definition (\ref{Ts}).
It then follows from the recursion (\ref{recur_taus}) that $\tau_{j}(1-x,n)=(-1)^j\tau_{j}(x,n)$ as well, and this implies that $|a_i(n,j)|=|a_{j-i}(n,j)|$ for all $i$.
Next, consider the polynomial $A(x):=\sum_{i=0}^j a_i(n,j)x^i$. Since $\tau_{j}(x,n)=(1-x)^j A(\frac{x}{1-x})$, Lemma \ref{lem_taus} implies that
$|A(x)| \le \const_j$ for $x \in [0,1/n]$. Thus $A_*(x):=A(\frac{x+1}{2n})$ satisfies $|A_*(x)| \le  \const_j$ for $x \in [-1,1]$.
Markov's inequality $\|A_*^{(i)}\|_\infty \le j^{2i} \|A_*\|_\infty$ (see \cite[Chapter 4, Theorem 1.4]{DeVoreL}) yields
$$ |a_i(n,j)|=(i!)^{-1} |A^{(i)}(0)|=(i!)^{-1}(2n)^i |A_*^{(i)}(-1)| \le C_j^\sharp n^i \mbox{ \rm for all } i \le j \,.
$$
\eop
For our next argument, we will need an additional technical lemma that provides bounds
on the derivatives of the functions $p_{nk}$.

\begin{lemma} \label{bounds_pnk}
For any integer $\ell \geq 0$ and any $\beta\geq 0$,
\begin{equation}
 \sum_{k=0}^n \left|  {k\over n} - x  \right|^\beta
\left| p^{(\ell)}_{nk}(x)   \right| \leq \const_{\beta,\ell} (\Delta_n(x))^{\beta-\ell}
\qquad \hbox{\rm for all \/} x\in [0,1].   \label{sum_deriv_pnk}
\end{equation}
\end{lemma}

\proof  The proof is by induction on $\ell$. The proof for $\ell=0$ is  due to
Lorentz \cite[Lemma 1]{Lorentz_paper}.
Our proof of the bound~(\ref{sum_deriv_pnk}) for $\ell \ge 1$ splits into
two cases. \vskip 2mm

\noindent
\underline{\bf Case 1}.  $\Delta_n(x)=\sqrt{x(1-x)/n}$. In this case
we start from the equality
$$ p'_{nk}(x)={k-nx \over x(1-x)} p_{nk}(x),  $$
and deduce by induction on $\ell$ that the $\ell$th derivative of $p_{nk}$ has the form
\begin{equation}
p_{nk}^{(\ell)}(x)=\sum_{i,j,\nu \ge 0} \Bigl\{ {n^i(k-nx)^\nu \over [x(1-x)]^j} \Psi_{\ell i j \nu} (x) p_{nk}(x)
\; : \;  i+j \le \ell ; \;  i+\nu \le j \Bigr\} \,, \label{pnk-ell}
\end{equation}
where $\Psi_{\ell i j \nu} (x)$ are polynomials in $x$ with coefficients that do not depend on $n$.
For fixed $i,j, \nu$, we have (using that $x(1-x) \ge \Delta_n(x)$ in this case)
$$ {n^i |k-nx|^\nu \over [x(1-x)]^j} =(\Delta_n(x))^{-2(i+\nu)} { | \frac{k}n-x|^\nu \over [x(1-x)]^{j-i-\nu}}
\le (\Delta_n(x))^{-(i+\nu+j)} \, \left|\frac{k}n-x\right|^\nu \, ,
$$
whence (using the already established case $\ell=0$ of (\ref{sum_deriv_pnk})), we  obtain
$$
\sum_{k=0}^n \left|  {k\over n} - x  \right|^\beta {n^i |k-nx|^\nu \over [x(1-x)]^j} p_{nk}(x) \le
(\Delta_n(x))^{-(i+\nu+j)} \sum_{k=0}^n \left|  {k\over n} - x  \right|^{\beta+\nu} p_{nk}(x)
\le \const_{\beta,\ell} (\Delta_n(x))^{\beta-i-j} \,.
$$
The restriction $i+j \le \ell$ implies that the right-hand side of the last display is at most
$\const_{\beta,\ell} (\Delta_n(x))^{\beta-\ell} \,.$ The representation (\ref{pnk-ell}) completes the proof in this case.

\vskip 2mm

\noindent
\underline{\bf Case 2}.  $\Delta_n(x)=1/n$. In this case we substitute a different expression for
$p'_{nk}(x)$, precisely
$$ p'_{nk}(x)=np_{n-1,k-1}(x)-np_{n-1,k}(x), $$
which yields
$$ p^{(\ell)}_{nk}(x)=np^{(\ell-1)}_{n-1,k-1}(x)-np^{(\ell-1)}_{n-1,k}(x) \, . $$
This gives the bound
$$ \sum_{k=0}^n \left|  {k\over n} - x  \right|^\beta
\left| p_{nk}^{(\ell)}(x)   \right|  \leq  \const_{\beta, \ell} \, n
\sum_{k=0}^n \left|  {k\over n} - x  \right|^\beta  \left(
\Big| p_{n-1,k-1}^{(\ell-1)}(x) \Big| + \Big|p_{n-1,k}^{(\ell-1)}(x) \Big|
\right). $$
The general inequality $(a+b)^\beta \le 2^\beta (a^\beta+b^\beta)$ implies that
\begin{eqnarray*}  \left|  {k\over n} - x  \right|^\beta \leq
2^\beta  \Bigl( \left|  {k-1\over n-1} - x  \right|^\beta
+\left(\frac{n-k}{n(n-1)}\right)^\beta \Bigr) \leq
2^\beta  \Bigl( \left|  {k-1\over n-1} - x  \right|^\beta +\frac{1}{n^\beta} \Bigr) \qquad {\rm for} \;\;\; k\geq 1,  \\
 \left|  {k\over n} - x  \right|^\beta \leq
2^\beta  \Bigl( \left|  {k\over n-1} - x  \right|^\beta +
\left(\frac{k}{n(n-1)}\right)^\beta \Bigr)  \leq 2^\beta  \Bigl( \left|  {k\over n-1} - x  \right|^\beta +\frac{1}{n^\beta} \Bigr) \qquad {\rm for} \;\;\; k<n.
\end{eqnarray*}
Therefore $\sum_{k=0}^n \Big|  {k\over n} - x  \Big|^\beta
\Big| p^{(\ell)}_{nk}(x)  \Big| $ is at most

\begin{eqnarray*}
2^\beta n\sum_{k=1}^n \Bigl( \left|  {k-1\over n-1} - x  \right|^\beta +n^{-\beta} \Bigr)
 | p_{n-1,k-1}^{(\ell-1)}(x)| + 2^\beta n \sum_{k=0}^{n-1} \Bigl( \Big|  {k\over n-1} - x  \Big|^\beta  +n^{-\beta} \Bigr) |p_{n-1,k}^{(\ell-1)}(x)|   \,.
\end{eqnarray*}
Since $\frac{n-1}n \Delta_{n-1}  \le \Delta_{n}  \le \Delta_{n-1}$,
we can finish the proof 
using the inductive assumption on $\ell-1$.

\eop

We now generalize Lemma~\ref{lem_taus} to derive bounds on the derivatives of the polynomials $\tau_{j}$.

\begin{lemma} \label{lem_derivs_taus} The derivatives of the polynomials $T_{nj}$ and
$\tau_{j}$ are bounded
as follows
\begin{eqnarray}
|T_{nj}^{(\ell)}(x)| & \leq  & \const_{j,\ell} \, n^j (\Delta_n(x))^{j-\ell} \label{derivs_T} \\
 | \tau^{(\ell)}_{j}(x,n)  | & \leq & \const_{j,\ell} \,  n^j (\Delta_n(x))^{j-\ell} \label{derivs_tau}
\end{eqnarray}
for all $x\in [0,1]$.
\end{lemma}

\proof Differentiate the formula~(\ref{Ts}) $\ell$ times to obtain
$$ T_{nj}^{(\ell)}(x) = \sum_{m\leq \min\{j, \ell\}}  {\ell \choose m}
\sum_{k=0}^n (-n)^m \frac{j!}{(j-m)!}(k-nx)^{j-m} p^{(\ell -m)}_{nk}(x).   $$
By Lemma~\ref{bounds_pnk}, each term is bounded by
$$ n^m \cdot n^{j-m} \cdot \const (\Delta_n(x))^{(j-m)-(\ell-m)}=
\const n^j (\Delta_n(x))^{j-\ell},$$
which proves the estimate~(\ref{derivs_T}).
To get the analogous estimate for derivatives of $\tau_{j}=\tau_{rj}$, we
run an inductive argument. Differentiating the formula~(\ref{recur_taus})
$\ell$ times, we get
\begin{equation}
 \tau^{(\ell)}_{j}(x,n)= - \sum_{m\leq \ell} {\ell \choose m} \sum_{s=2}^j
{1\over s!} T^{(m)}_{ns}(x) \tau^{(\ell - m)}_{j-s}(x,n).
\label{sum_der_taus}
\end{equation}
Applying the inductive assumption on the derivatives $\tau^{(\ell - m)}_{j-s}(x,n)$
and the already proven bound~(\ref{derivs_T}) on $T^{(m)}_{ns}(x)$, we obtain the
estimate
$$  \const_{j,\ell} n^s (\Delta_n(x))^{s-m}\, n^{j-s} (\Delta_n(x))^{j-s-\ell+m}=
\const_{j,\ell} n^j (\Delta_n(x))^{j-\ell}  $$
on each term in the sum~(\ref{sum_der_taus}), and therefore on the function
$|\tau^{(\ell)}_{j}(x,n)|$ as well, proving~(\ref{derivs_tau}).
\eop

Next, we will show that the derivatives of the polynomials $Q_{n,r}f$ approximate the
corresponding derivatives of $f$ sufficiently well. This is known as {\bf\em simultaneous
approximation}. Here is the precise result.

\begin{lemma}  \label{lem-4}
Let $f\in C^\alpha[0,1]$ and let $r\eqbd \aint$.  Then, for any $j=0, \ldots, r$,
$$ |((I-Q_{n,r})f)^{(j)}(x)|\leq C_{\ref{lem-4}} \|f \|_{C^\alpha} (\Delta_n(x))^{\alpha-j}
\qquad \hbox{\rm for all} \;\; x\in [0,1],   $$
where the constant $C_{\ref{lem-4}}$ is independent of $f$ and $n$.
\end{lemma}

\proof The case $j =0$ of this lemma is contained in formula (22) of \cite{Lorentz_paper}.
To prove the result  for all $j$, we use the fact the Lorentz
operator $Q_{n,r}$ reproduces polynomials of degree at most $r$.

Our goal is to show that the $j$th derivative of the difference between $Q_{n,r}f$
and $f$ at any point $x$ is bounded by a constant multiple of $\|f \|_{C^\alpha}
(\Delta_n(x))^{\alpha-j}$ regardless of $x$. Since $Q_{n,r}$ reproduces polynomials
of degree $r$ (by Lemma~\ref{lem_rep}), we can subtract from $f$ its Taylor polynomial
of degree $r$ centered at $x$ without changing the difference $(Q_{n,r}f-f)^{(j)}(x)$.
Thus, without loss of generality we can assume that the value of $f$ and its derivatives
up to order $r$ are zero at $x$. Now, recall that
$$ (Q_{n,r}f-f)^{(j)}(x) = \left( \sum_{k=0}^n \left(f\left({k\over n}\right) -f(x)\right)
 p_{nk}(x) +\sum_{i=2}^r \sum_{k=0}^n  f^{(i)}\left({k\over n}\right) {1\over n^i}
\tau_{i}(x,n) p_{nk}(x) \right)^{(j)}.   $$
Differentiating these sums $j$ times, we will obtain terms of two
kinds. Terms of the first kind are obtained from differentiating the
first sum; they have the form
$$ \sum_{k=0}^n \left(  f\left({k\over n}\right) -f(x)  \right)^{(\ell)} p_{nk}^{(j-\ell)}(x)
$$ for some $\ell$ between $0$ and $j$. Each of these sums can be bounded
as follows, using Lemma \ref{bounds_pnk}:
\begin{eqnarray*}
 \left| \sum_{k=0}^n \left( f\left({k\over n}\right)-f(x)  \right)^{(\ell)}
p_{nk}^{(j-\ell)}(x)   \right| & \leq & \sum_{k=0}^n \|f \|_{C^\alpha}
\left|{k\over n}-x\right|^{\alpha-\ell}  |p_{nk}^{(j-\ell)}(x) | \\
& \leq & \const \| f \|_{C^\alpha} (\Delta_n(x))^{(\alpha-\ell ) -(j-\ell )} \\ & = &
\const \|f \|_{C^\alpha}  (\Delta_n(x))^{\alpha-j} .
\end{eqnarray*}
Terms of the second kind are obtained by differentiating any of the other
sums for $i=2, \ldots, r$ and have the form
$$  \sum_{k=0}^n f^{(i)}\left({k\over n} \right) {1\over n^i}
\tau_{i}^{(\ell)} (x,n) p_{nk}^{(j-\ell)}(x)  $$ for some $\ell$ between
$0$ and $j$. Taking into account that the derivatives of $f$ up
to order $r$ vanish at $x$, each of these sums can be bounded by
\begin{eqnarray*}
\left| \sum_{k=0}^n f^{(i)}\left({k\over n} \right) {1\over n^i}
\tau_{i}^{(\ell)} (x,n) p_{nk}^{(j-\ell)}(x) \right| \leq
\sum_{k=0}^n \|f \|_{C^\alpha}
 \left|  {k\over n} - x   \right|^{\alpha - i} {1\over n^i}
 |\tau_{i}^{(\ell)}(x,n)| \, \left| p_{nk}^{(j-\ell)}(x)  \right|.
\end{eqnarray*}
Invoking the bound~(\ref{derivs_tau}) from Lemma~\ref{lem_derivs_taus} on
the terms $|\tau_{i}^{(l)}(x,n) |$, we conclude that the total is bounded by
$$ \const \|f \|_{C^\alpha} {1\over n^i} \, n^i (\Delta_n(x))^{i-\ell}
\sum_{k=0}^n \left| {k\over n} - x \right|^{\alpha-i} |p^{(j-\ell)}_{nk}(x)|. $$
The last sum, in turn, is estimated according to Lemma~\ref{bounds_pnk} to
produce the final bound
$$ \const \| f\|_{C^\alpha}
(\Delta_n(x))^{i-\ell} \cdot \Delta^{\alpha-i-j+\ell}_n(x) =
\const \|f \|_{C^\alpha} (\Delta_n(x))^{\alpha-j}. $$
This completes the proof.
\eop

\begin{lemma}\label{lem-3} Let $f\in C^\alpha[0,1]$ and let
 $r\eqbd \lceil \alpha \rceil -1$. Then, for all $x \in [0,1]$,
\begin{equation}
 |  \left(Q_{n,r}f \right)^{(r+1)}(x) |   \leq   C_{\ref{lem-3}} 
(\Delta_n(x))^{\alpha-r-1}  \|f \|_{C^\alpha}~, \label{C-r+1-bound}
\end{equation}
with the  constant $C_{\ref{lem-3}}$ independent of $f$ and $n$.
\end{lemma}

\proof    Firstly, we may assume without loss of
generality that $f$ vanishes to order $r$ at $x$, since polynomials
of degree at most $r$ are reproduced by $Q_{n,r}$ (Lemma~\ref{lem_rep}) and then annihilated
by taking the derivative of order $r+1$, as well as by taking the $r$th derivative followed by a difference at two points
$x$ and $y$. The assumption made above implies that, for all $i \le r$,
$$ \Bigl|f^{(i)} \Bigl({k \over n}\Bigr)\Bigr| \leq \const \|f \|_{C^\alpha} \, \Bigl|\frac{k}n-x\Bigr|^{\alpha-i}. $$
By direct differentiation of (\ref{simpleQs}),
\begin{equation} \label{Qexp}(Q_{n,r}f)^{(r+1)}(x) =
\sum_{i=0}^r \sum_{\ell=0}^{r+1} {r+1 \choose \ell} \sum_{k=0}^n  f^{(i)}\left({k\over n}\right) {1\over n^i}
\tau_{i}^{(\ell)}(x,n) p_{nk}^{(r+1-\ell)}(x) .
\end{equation}
Fix $i \in [0,r]$ and $\ell \in [0,r+1]$. The summand corresponding to $i$ and $\ell$ in (\ref{Qexp}) can be bounded by
\begin{eqnarray} \label{tot}
\left| \sum_{k=0}^n f^{(i)}\left({k\over n} \right) {1\over n^i}
\tau_{i}^{(\ell)} (x,n) p_{nk}^{(r+1-\ell)}(x) \right| \leq
\sum_{k=0}^n \|f \|_{C^\alpha}
 \left|  {k\over n} - x   \right|^{\alpha - i} {1\over n^i}
 |\tau_{i}^{(\ell)}(x,n)| \, \left| p_{nk}^{(r+1-\ell)}(x)  \right|.
\end{eqnarray}
Invoking Lemma~\ref{lem_taus}, we note that the terms ${1\over n^i}|\tau_{i}^{(\ell)}(x,n) |$
are bounded by a constant multiple of $(\Delta_{n}(x))^{i-\ell}$, therefore
(\ref{tot}) is bounded by
$$ \const \|f \|_{C^\alpha}  (\Delta_{n}(x))^{i-\ell} (\Delta_n(x))^{\alpha-i-r-1+\ell}=
\const \|f \|_{C^\alpha} (\Delta_n(x)^){\alpha-r-1}. $$
This proves (\ref{C-r+1-bound}).   \eop

\begin{lemma}\label{lem-2}  Let $f\in C^\alpha[0,1]$ and let
 $r\eqbd \lceil \alpha \rceil -1$. Then, for any $x \in [0,1]$,
\begin{equation}
\|Q_{n,r} f\|_{C^{\alpha}} \leq C_{\ref{lem-2}} \|f\|_{C^{\alpha}}~,   \label{C-alpha-bound}
\end{equation}
with the  constant $C_{\ref{lem-2}}$ independent of $f$ and $n$.
\end{lemma}

\proof  To establish the  bound~(\ref{C-alpha-bound}), we need to estimate
the expression
\begin{equation}
 \left|  (Q_{n,r}f)^{(r)}(x) - (Q_{n,r}f)^{(r)}(y) \right|
\label{Qdiff}
\end{equation}
for two points $x$ and $y$ in $[0,1]$.  With loss of generality,
we may assume that $\Delta_n(x) \ge \Delta_n(y)$. Consider two cases.

\noindent\underline{\bf Case 1}. If  $|x-y| \geq \Delta_n(x)$, then we  estimate~(\ref{Qdiff}) using the triangle inequality and the bound
$$ | (Q_{n,r} f - f)^{(r)}(x)| \leq \const (\Delta_n(x))^{\alpha-r} \|f \|_{C^\alpha}$$
from Lemma \ref{lem-4} on each of the two terms,  $(Q_{n}f-f)^{(r)}(x)$ and $(Q_{n}f-f)^{(r)}(y)$.
Altogether, this bounds (\ref{Qdiff}) from above by
$$ \const \|f \|_{C^\alpha} (\Delta_n(x))^{r-\alpha} \leq
 \const \|f \|_{C^\alpha} |x-y|^{r-\alpha}. $$

\noindent\underline{\bf Case 2}.  If $|x-y|\leq \Delta_n(x)$,  then  $|x-y| \le 
(\Delta_n(x))^{r+1-\alpha}  |x-y|^{\alpha-r} $ so for some $\xi$ between $x$ and $y$,
\begin{eqnarray} \nonumber
 \left|  (Q_{n,r}f)^{(r)}(x) - (Q_{n,r}f)^{(r)}(y) \right| &=& (Q_{n,r}f)^{(r+1)}(\xi) \cdot |x-y| \\
 &\le &
  \const (\Delta_n(\xi))^{\alpha-r-1} \|f \|_{C^\alpha} \cdot (\Delta_n(x))^{r+1-\alpha} |x-y|^{\alpha-r} \,.
\label{Qdiff2}
\end{eqnarray}
 Lemma \ref{lem-1} implies that $\Delta_n(x) \le 4\Delta_n(\xi)$, and inserting this bound in (\ref{Qdiff2})
 establishes (\ref{C-alpha-bound}).   \eop


\begin{lemma}  \label{lem-7} Suppose that  $f:[0,1] \to \R$
satisfies
 $|f^{(r+1)}(x)| \le (\Delta_n(x))^{-\beta}$ for some $\beta \in [0,1]$ and all $x \in [0,1]$. Then,  
for all $x \in [0,1]$, we have
\begin{equation}
   \left|  \Bigl(Q_{n,r} f \Bigr)^{(r+1)}(x) \right| \leq    C_{\ref{lem-7}}
 (\Delta_n(x))^{-\beta} \,.
\label{hybrid}
\end{equation}
with  $C_{\ref{lem-7}}=C_{\ref{lem-7}}(r,\beta)$ a constant independent of $f$ and $n$.
\end{lemma}

\proof
To prove (\ref{hybrid}), we may assume as in the preceding theorem
that $f$ vanishes to order $r$ at $x$.
This  implies that for all $i \le r$ and $z \ne x$ in $[0,1]$, there exists $\xi$ between $x$ and $z$ such that
\begin{equation} \label{tayl2} \frac{|f^{(i)}(z)|}{|z-x|^{r+1-i}} \leq  |f^{(r+1)}(\xi)| \le 
(\Delta_n(\xi))^{-\beta} \le 2(\Delta_n(x))^{-\beta} \Bigl(1+\frac{|x-z|}{\Delta_n(x)}\Bigr)\,,
\end{equation}
where the last step used (\ref{Delta-bound}) taken to the power $\beta$, and the inequality 
$|x-\xi| \le |x-z|$.

Recall the expression (\ref{Qexp}) for $(Q_{n,r}f)^{(r+1)}(x)$.
Fix $i \in [0,r]$ and $\ell \in [0,r+1]$. The summand
$$
\left| \sum_{k=0}^n f^{(i)}\left({k\over n} \right) {1\over n^i}
\tau_{i}^{(\ell)} (x,n) p_{nk}^{(r+1-\ell)}(x) \right|
$$
corresponding to $i$ and $\ell$ in (\ref{Qexp}) can be bounded using  (\ref{tayl2}) and Lemma~\ref{lem_taus} by
\begin{eqnarray} \label{tot2}
\const \sum_{k=0}^n
 \left|  {k\over n} - x   \right|^{r+1 - i}  (\Delta_n(x))^{-\beta} \Bigl(1+\frac{|{k \over n}-x|}
{\Delta_n(x)}\Bigr)\,(\Delta_{n}(x))^{i-\ell} \,  \left| p_{nk}^{(r+1-\ell)}(x)  \right|.
\end{eqnarray}
Invoking Lemma \ref{bounds_pnk} twice, we conclude that
(\ref{tot2}) is bounded by
$$ \const   (\Delta_{n}(x))^{-\beta} (\Delta_n(x))^{i-\ell} \Bigl((\Delta_n(x))^{\ell-i}+
\frac{(\Delta_n(x))^{\ell-i+1}}{\Delta_n(x)}\Bigr) \,
\le \const   (\Delta_{n}(x))^{-\beta} \,.
$$
This proves the lemma.
\eop

\section{The iterative construction} \label{sec-iterate}

The goal of this section is to prove the suffiency part of Theorem~\ref{thm-main}.
This will be achieved via an iterative construction of the approximants $f_n$
that are subsequently adjusted downward and upward to produce the approximants
$g_n$ and $h_n$ satisfying the required consistency conditions listed in
Result~\ref{res_reduction} in the Introduction. We begin by  analyzing the
behaviour of the degree $n+r$ Bernstein coefficients of $Q_{n,r} f$.

\begin{lemma}\label{lem-5}
For every $\epsilon > 0$, there exists $n_0$ such that for $n \ge n_0$, the
degree $n+r$ Bernstein coefficients of $Q_{n,r} f$ are between $\min_{[0,1]} f 
- \epsilon$ and $\max_{[0,1]} f + \epsilon$~.
\end{lemma}

\proof Recall that
$$Q_{n,r} f = \sum_{k=0}^n \left(\sum_{j=0}^r \frac{f^{(j)}(\frac{k}{n})}{n^j} \tau_j(x,n)\right)p_{n,k}(x)~.$$
Note that the $i$th Bernstein coefficient of $\tau_j(x,n)/n^j \in {\cal B}_j$ is bounded
by $C_j^\sharp \min(n^{-i}, n^{i-j})$ by Corollary~\ref{cor_taus}. This  implies
that the  Bernstein coefficients of  $f^{(j)}(\frac{k}{n}) \tau_j(x,n)/n^j$ for
$j\geq 1$ do not exceed $\frac{\const}{n}\max_{1 \le j \le r} \|f^{(j)}\|_{\infty}$. Since the
Bernstein coefficients of $\sum_{k=0}^n f(\frac{k}{n})p_{n,k}(x)$
are between $\min_{[0,1]} f$ and $\max_{[0,1]} f$, we conclude that
the Bernstein coefficients of $Q_{n,r} f \in {\cal B}_{n+r}$ are between
$\min_{[0,1]} f - \frac{\const}{n} \max_{1 \le j \le r} \|f^{(j)}\|_{\infty}$ and $\max_{[0,1]} f +
\frac{\const}{n} \max_{1 \le j \le r} \|f^{(j)}\|_{\infty}$, from which Lemma \ref{lem-5} follows
immediately.  \eop

\begin{lemma}\label{lem-6} Let $r\eqbd \lceil \alpha \rceil -1$.
If $f^{(j)}(x) \leq (\Delta_n(x))^{\alpha-j}$ $(j=0, \ldots, r)$ for all $x\in [0,1]$, 
then the degree $n+r$ Bernstein coefficients of $Q_{n,r} f$ are dominated by those of
$$C_{\ref{lem-6}} [x+ (1-x)]^r \left[\sum_{k=0}^n \left(\Delta_n\left(\frac{k}{n}\right)\right)^\alpha 
p_{n,k}(x)\right]~,$$
where $C_{\ref{lem-6}}$ does not depend on $n$ and $f$.
\end{lemma}

\proof This lemma is a bit trickier. Separating the contributions given by
different $j = 0, \ldots, r$, we see that it would suffice to bound
the coefficients of
$$\left[\sum_{k=0}^n \left(\Delta_n\left(\frac{k}{n}\right)\right)^{\alpha - j}  p_{n,k}(x)\right] 
\cdot {\tau_j (x,n) \over n^j} \in {\cal B}_{n+j}$$
by those of $\big[\sum_{k=0}^n \left(\Delta_n\left(\frac{k}{n}\right)\right)^\alpha 
p_{n,k}(x)][x+(1-x)\big]^j$ (possibly with some constant
factor). Since both polynomials have symmetric coefficients as
Bernstein polynomials in ${\cal B}_{n+r}$ and since we may assume without
loss of generality that $n>3r$, we see that it is enough to prove
that
\begin{align*}
&\sum_{\substack{s+t=u \\ 0\leq s\leq u \\ 0\leq t\leq j}}{n\choose s}
\left(\Delta_n\left(\frac{s}{n}\right)\right)^{\alpha - j}
\min\{ n^{-t}, n^{t-j} \}  \leq  C \sum_{\substack{s+t =u \\ 
0\leq s\leq u\\ 0\leq t\leq j}} {n\choose s}
\left( \Delta_n \left(\frac{s}{n}\right)\right)^\alpha~ 
\quad \mbox{for }\; 0\leq u\leq \frac{n+r}{2} \leq \frac{2n}{3}~.
\end{align*}
Note  that $\Delta_n\left(\frac{s}{n}\right)$ is comparable
to $\Delta_n\left(\frac{u}{n}\right)$ for $|s-u| \leq j\leq r$. This
allows us to reduce the inequality to
$$\sum_{\substack{s+t=u \\ 0\leq s \leq u \\ 0\leq t\leq j}} {n \choose s} \min\{
n^{-t}, n^{t-j} \} \leq \const \Big[\sum_{\substack{s+t = u\\ 0\leq s\leq u\\
 0\leq t\leq j}}{n\choose s}\Big] \Delta_n^j\left(\frac{u}{n}\right)~.$$
We shall keep just one term ${n \choose u} \left(\Delta_n\left(\frac{u}{n}\right)\right)^j$ 
on the right and use the estimate ${n\choose s} = {n \choose u-t} \leq 3^t
\left(\frac{u}{n}\right)^t{n \choose u}$ valid for $u \leq \frac{2n}{3}$, $t\geq 0$. 
Since $t\leq j \leq r$, we only need to show that
\begin{align*}
\left(\frac{u}{n}\right)^t \min\{ n^{-t}, n^{t-j}\} =
\min\left\{ \frac{u^t}{n^{2t}}, \frac{u^t}{n^j}\right\} \leq C 
\left(\Delta_n\left(\frac{u}{n}\right)\right)^j~.
\end{align*}
But $\left(\Delta_n\left(\frac{u}{n}\right)\right)^j \geq \frac{1}{n^j}$, which takes
care of $u=0$ (with the only possible $t=0$), and if $1\leq u \leq \frac{2n}{3}$, we 
have $\Delta_n\left(\frac{u}{n}\right) \geq C\frac{u^{1/2}}{n}$, so it suffices to prove 
that $\min \left\{ \frac{u^t}{n^{2t}}, \frac{u^t}{n^j} \right\} \leq
\frac{u^{j/2}}{n^j}$ or equivalently,
$\min\{\left(\frac{u}{n^2}\right)^{t-\frac{j}2}, u^{t-\frac{j}2}\} \leq 1$. But the
first term is less than 1 for $t>\frac{j}{2}$ and the second term is
not greater than 1 for $t\leq \frac{j}{2}$.  \eop

\subsection*{Iterative construction of $f_n$}

 Let $\alpha > 0$, $\alpha \not\in \mathbb{Z}$, hence $r= \lfloor\alpha\rfloor$. Assume that $f\in C^\alpha [0,1]$ satisfies
$$0< \min_{[0,1]} f \leq \max_{[0,1]} f < 1 \,.
$$
 Fix $n_0 \in \N$ and  $b = 2^s$ to be chosen later. Denote $\Lambda \eqbd \{b^m n_0 : m\geq 0\}$ and define $f_n$ for $n\in \Lambda$ by
 \begin{eqnarray*}
f_{n_0} & \eqbd & Q_{n_0,r}f~,  \\
f_n & \eqbd & f_{n/b} + Q_{n,r}(f- f_{n/b})\quad  \hbox{ for }\; n > n_0~.
\end{eqnarray*}
Our task is to show that $f_n \rightarrow f$, that the Bernstein coefficients of $f_n$ are between $\delta$ and $1-\delta$ for some $\delta > 0$, and that the Bernstein coefficients of $Q_{n,r}(f - f_{n/b})$ are dominated (up to some constant factor) by those of
 $$[x + (1-x)]^r \Big[\sum_{0\leq k \leq n} \left(\Delta_n \left(\frac{k}{n}\right)\right)^\alpha
p_{n,k}(x)\Big]~.$$
We will do it in four steps.  \medskip

\noindent
\underline{\textbf{Step 1}}. {\bf Estimate for $f_n^{(r+1)}$.}  \/ We will show by induction that
 \begin{equation} \label{astep1}
 |f_n^{(r+1)}| \leq 2 C_{\ref{lem-3}} \|f\|_{C^\alpha} (\Delta_n)^{\alpha -r -1} \; 
\mbox{ \rm on } \; [0,1],  \end{equation}
 provided that $b$ is chosen large enough. By Lemma \ref{lem-3}, the inequality (\ref{astep1}) holds for $n=n_0$. Assume that it is true for $n/b$ in place of $n$. Write
  $$|f_n^{(r+1)}| \leq |f_{n/b}^{(r+1)}| + |(Q_{n,r} f_{n/b})^{(r+1)}| + |(Q_{n,r} f)^{(r+1)}|~.$$
According to Lemma \ref{lem-3}, the last term is bounded by
 $C_{\ref{lem-3}} \|f\|_{C^\alpha} (\Delta_n)^{\alpha-r-1}$. By the induction hypothesis,
$$|f_{n/b}^{(r+1)}| \leq 2 C_{\ref{lem-3}} \|f\|_{C^\alpha} (\Delta_{n/b})^{\alpha - r -1} 
\leq 2 b^{(\alpha - r -1)/2} C_{\ref{lem-3}} \|f\|_{C^\alpha} (\Delta_n)^{\alpha - r -1}$$
whence by Lemma \ref{lem-7} (with proper renormalization)
$$|(Q_{n,r} f_{n/b})^{(r+1)}| \leq 2 C_{\ref{lem-7}} b^{(\alpha - r -1)/2} C_{\ref{lem-3}} \|f\|_{C^\alpha} (\Delta_n)^{\alpha - r -1}~.$$
If $b$ is chosen so large that $2(1+C_{\ref{lem-7}})b^{(\alpha - r -1)/2} \leq 1$, we can add these three estimates to get
$$|f_n^{(r+1)}| \leq 2 C_{\ref{lem-3}} \|f\|_{C^\alpha} (\Delta_n)^{\alpha -r-1}~.$$
Moreover, we see that

\begin{equation} \label{zstep1}
|[(I-Q_{n,r})f_{n/b}]^{(r+1)}| \leq 2(1+C_{\ref{lem-7}}) b^{(\alpha - r- 1)/2}C_{\ref{lem-3}} \|f\|_{C^\alpha} (\Delta_n)^{\alpha - r -1}~.
\end{equation}  \medskip

\noindent
\underline{\bf{Step 2}}. {\bf An estimate for $\|f_n\|_{C^\alpha}$.} \/
We will show that
$$\|f_n\|_{C^\alpha} \leq 2 C_{\ref{lem-2}} \|f\|_{C^\alpha}, \; \mbox{ \rm provided  that $b$ is large enough.} $$
 Again, we will argue by induction. Lemma \ref{lem-2} yields the base case $n = n_0$. Assume that the statement is true for $n/b$. Write
$$f_n = (I- Q_{n,r})f_{n/b} + Q_{n,r} f~.$$
Since  $\|Q_{n,r} f\|_{C^\alpha} \leq C_{\ref{lem-2}} \|f\|_{C^\alpha}$ by Lemma \ref{lem-2}, it suffices to show that the $C^\alpha$-norm of the function 
$$\Psi \eqbd (I - Q_{n,r}) f_{n/b}$$ is bounded by $C_{\ref{lem-2}}\|f\|_{C^\alpha}$.

We need to estimate $|\Psi^{(r)}(x) - \Psi^{(r)}(y)|$. Without loss of generality, we may assume that $\Delta_n(x) \geq \Delta_n (y)$. Choose a big positive constant $A$ and consider two cases: \smallskip

\noindent{\bf Case 1.} $|x-y| \geq A \Delta_n(x)$. Then
\begin{eqnarray*}
|\Psi^{(r)} (x) - \Psi^{(r)} (y)| & \leq & |\Psi^{(r)}(x)| + |\Psi^{(r)}(y)| \\
 &\leq & 2C_{\ref{lem-4}} C_{\ref{lem-2}} \|f\|_{C^\alpha} ((\Delta_n(x))^{\alpha - r} 
+ (\Delta_n(y))^{\alpha -r})  \\
& \leq & 4  C_{\ref{lem-4}} C_{\ref{lem-2}} \|f\|_{C^\alpha} (\Delta_n(x))^{\alpha - r}  \\
&\leq & 4C_{\ref{lem-4}} A^{-(\alpha - r)} C_{\ref{lem-2}} \|f\|_{C^\alpha} |x-y|^{\alpha - r}
\end{eqnarray*}
and we get the desired bound if $ C_{\ref{lem-4}} A^{-(\alpha - r)} \leq 1$.
\smallskip

\noindent{\bf Case 2.} $|x - y| \leq A \Delta_n(x)$.
For this case, we will use the estimate
$$|\Psi^{(r +1)}| \leq 2(1+C_{\ref{lem-7}})b^{(\alpha - r -1)/2} C_{\ref{lem-3}} 
\|f\|_{C^\alpha} (\Delta_n)^{\alpha - r-1}$$
obtained in (\ref{zstep1}). Write
$$|\Psi^{(r)} (x) - \Psi^{(r)} (y)| =  |\Psi^{(r+1)} (\xi)| |x-y|$$
for some $\xi$ between $x$ and $y$.
Now, by Lemma \ref{lem-1}, $\Delta_n (\xi) \geq [2(1+A)]^{-1} \Delta_n(x)$. Combining this with the above estimate for $|\Psi^{(r+1)}|$, we obtain
\begin{eqnarray*}
|\Psi^{(r+1)} (\xi)| |x-y| & \leq & 2(1+C_{\ref{lem-7}}) b^{(\alpha  -r -1)/2}C_{\ref{lem-3}} 
\|f\|_{C^\alpha} [2(1+A)]^{r+1 -\alpha}A^{r+1-\alpha}[A\Delta_n(x)]^{\alpha - r -1}|x - y|\\
& \leq& 2 (1+C_{\ref{lem-7}}) [2A(1+A)]^{r+1 - \alpha} b^{(\alpha -r-1)/2}C_{\ref{lem-3}} 
\|f\|_{C^\alpha} |x-y|^{\alpha -r}
\end{eqnarray*}
and we get the desired conclusion if
$$2(1+ C_{\ref{lem-7}}) [2A(1+A)]^{r+1 -\alpha} b^{(\alpha - r - 1)/2} \leq 1~.$$
\smallskip

\noindent
\underline{\bf{Step 3}}. {\bf An estimate for $(f-f_n)^{(j)}$}.  \/
Since $f - f_n = (I -Q_{n,r}) (f - f_{n/b})$ for $n \geq bn_0$ and we know that 
$\|f_{n/b}\|_{C^\alpha} \leq 2 C_{\ref{lem-2}} \|f\|_{C^\alpha}$, we can invoke 
Lemma \ref{lem-4} to conclude that
$$|(f - f_n)^{(j)}| \leq C_{\ref{lem-4}} (1+ 2 C_{\ref{lem-2}})\|f\|_{C^\alpha} 
(\Delta_n)^{\alpha - j} \quad \hbox{  for  } \;\; n \geq b n_0~.$$
The same, or an even better, estimate can be derived for $n = n_0$ from the 
representation $f-f_{n_0} = (I - Q_{n_0})f$. In particular, we see that $f_n \to f$ 
uniformly in $[0,1]$.
\medskip

\noindent
\underline{\bf{Step 4}}. {\bf Estimates for  Bernstein coefficients.}  \/
It follows now from Lemma \ref{lem-6} and the result of the previous step that the degree $n+r$
Bernstein coefficients of $Q_{n,r} (f -f_{n/b})$ are dominated by those of
$$C_{\ref{lem-6}} C_{\ref{lem-4}} (1 + 2 C_{\ref{lem-2}}) b^{\alpha /2} \|f\|_{C^\alpha} [x+(1-x)]^r \Big[\sum_{k=0}^n (\Delta_n)^\alpha \left(\frac{k}{n}\right) p_{n,k}(x)\Big]$$
(here we used the inequality $(\Delta_{n/b})^j\leq b^{j/2} (\Delta_n)^j \leq b^{\alpha /2} (\Delta_n)^j$ 
for $ 0 \leq j  \leq r$).
Since the latter coefficients are bounded by
$$C_{\ref{lem-6}} C_{\ref{lem-4}} (1+2C_{\ref{lem-2}}) b^{\alpha /2} \|f\|_{C^\alpha} n^{-\alpha /2}~,$$
we see that the degree $n+r$ Bernstein coefficients of $f_n$ differ from those of $f_{n_0}$ at most by 
the factor $$\const \sum_{n\in \Lambda, n > n_0} n^{-\alpha/2} \leq C_* \cdot n_0^{-\alpha/2}.$$
 Now, fix $\delta > 0$ such that
$$
2\delta < \min_{[0,1]} f \leq \max_{[0,1]} f < 1- 2\delta
$$ and choose $n_0$  large enough so that $C_* \cdot n_0^{-\alpha/2} < \delta$ and the degree $n+r$ 
coefficients of $Q_{n_0} f$ are between $2 \delta$ and $(1-2\delta)$, which is possible by 
Lemma \ref{lem-5}. Then the degree $n+r$ Bernstein coefficients of $f_n$ are between $\delta$ and 
$1-\delta$ for all $n \in \Lambda$ such that $n \ge n_0$.
\medskip

\noindent \underline{\bf Step 5}.   {\bf{Construction of $g_n$ and $h_n$}}.
Set
\begin{equation}
 \varphi_n(x) \eqbd  \frac{\theta_\alpha}{n^\alpha} + \left[\frac{x(1-x)}{n}\right]^{\alpha/2} \,,
\label{phi}
\end{equation} where $\theta_\alpha$ will be specified later, and define
\begin{equation} \label{defgh}
g_n  \eqbd f_n - [x + (1-x)]^r DB_n\varphi_n~, \quad h_n \eqbd f_n + [x+(1-x)]^r DB_n\varphi_n~.
\end{equation}
The constant $D$ here is to be chosen later.
Clearly, the degree $n{+}r$ Bernstein coefficients of $h_n$ are greater than those of $g_n$. Also, since $|\varphi_n| < \delta D^{-1}$ for large $n$, we see that the Bernstein coefficients of $g_n$ are positive and those of $h_n$ are less than 1 for sufficiently large $n$. It remains to show that the Bernstein coefficients of $g_n$ ``increase'', those of $h_n$ ``decrease'' and that $g_n - h_n = O ((\Delta_n)^\alpha)$.

\begin{lemma}\label{lem-X}  The functions $\varphi_n$ defined in~(\ref{phi}) satisfy
\begin{equation*}
B_n \varphi_n \leq \const \varphi_n \leq \const (\Delta_n)^\alpha~.
\end{equation*}
Consequently,  $h_n - g_n = O ((\Delta_n)^\alpha)$.
\end{lemma}

\proof
Since $\varphi_n$ is comparable to $(\Delta_n)^{\alpha}$, it
suffices to show that
$$
\sum_{k=0}^{n} \left(\Delta_n \left(\frac{k}{n}\right)\right)^\alpha p_{n,k}(x) \leq C 
(\Delta_n(x))^\alpha~.
$$
 Recall that, by Lemma \ref{lem-1},
$$\left(\Delta_n\left(\frac{k}{n}\right)\right)^\alpha \leq 2^\alpha 
\left(1+ \frac{\left|x - \frac{k}{n}\right|}{\Delta_n(x)}\right)^\alpha 
(\Delta_n(x))^\alpha \leq 2^{2\alpha} \left((\Delta_n(x))^\alpha + 
\left|x - \frac{k}{n}\right|^\alpha\right)~.$$
Now, the first term yields the sum
$$2^{2\alpha} (\Delta_n(x))^\alpha \sum_{k=0}^n p_{n,k}(x) = 2^{2\alpha} (\Delta_n(x))^\alpha~,$$
while the second one yields the sum
$$2^{2\alpha} \sum_{k=0}^n \left|x - \frac{k}{n}\right|^\alpha p_{n,k}(x) \leq \const 
(\Delta_n(x))^\alpha \,, $$ due to Lemma~\ref{bounds_pnk}.  The desired bound
$h_n - g_n = O ((\Delta_n)^\alpha)$ now follows from (\ref{defgh}).
\eop

Next, we want to show that
\begin{equation} \label{want}
(B_n \varphi_n)  \succeq_{2n} (1+\gamma) B_{2n} \varphi_{2n} \,,
\end{equation}
with some $\gamma > 0$.  To perform this comparison, we multiply the
left-hand side by $[x+(1-x)]^n$ and expand. We thus see that this
claim is equivalent to the system of inequalities
$$\sum_{j=0}^k \frac{{n \choose j}{n \choose k-j}}{{2n \choose k}} \varphi_n\left(\frac{j}{n}\right) \geq (1 + \gamma )\varphi_{2n} \left(\frac{k}{2n}\right),~~\hbox{ for } 0\leq k \leq 2n~.$$
Denote the coefficients $\frac{{n\choose j} {n\choose k-j}}{{2n
\choose k}}$ by $\sigma_{k,j}$. Note that $\varphi_n \geq
2^{\alpha/2} \varphi_{2n}$, which immediately takes care of $k=0$
and $k=2n$ with any $\gamma < 2^{\alpha /2} -1$. So, we will assume
below that $1 \leq k \leq 2n-1$.

Note that the function $$ \Upsilon(x) \eqbd [x(1-x)]^{\alpha/2} $$ satisfies the
inequality
\begin{equation} \label{upsin}
\frac{\Upsilon(x+t) + \Upsilon(x-t)}{2} \geq \Upsilon(x) \left[1- \frac{c_\alpha}{\min\{ x, 1-x\}^2} t^2\right]~,
\end{equation}
for $0\leq t\leq \min\{x, 1-x\}$, provided that $c_\alpha$ is large enough.

Indeed, when $0\leq t\leq \frac{1}{2} \min\{ x, 1-x\}$, (\ref{upsin})  follows from
the estimate $|\Upsilon''(\xi)| \leq \const \frac{\Upsilon(x)}{\min\{x, 1-x\}^2}$, valid for all $\xi\in[x-t, x+t]$,
and when $$ \frac{1}{2} \, \min\{ x, 1-x\} <t \le \min\{ x, 1-x \} \,,$$
the inequality (\ref{upsin}) is trivial, provided that $c_\alpha \geq
4$. Taking into account that $\sigma_{k,j} = \sigma_{k,k-j}$,
$\sum_j \sigma_{k,j} =1$ and that $\sigma_{k,j} =0$ if $|\frac{j}{n}
-\frac{k}{2n}| > \min\{ \frac{k}{2n}, 1-\frac{k}{2n}\}$, we obtain
\begin{align*}
\sum_j \sigma_{k,j}\Upsilon\left(\frac{j}{n}\right) &\geq \Upsilon\left(\frac{k}{2n}\right)
\left[1- \frac{c_\alpha n^2}{\min\{k, 2n-k\}^2} \sum_j \sigma_{k,j}\left(\frac{j}{
n} - \frac{k}{2n}\right)^2\right]\\
&= \Upsilon\left(\frac{k}{2n}\right) \Big[1- \frac{c_\alpha k(2n-k)}{\min\{ k , 2n-k\}^2
4(2n-1)}\Big] \\
&\geq \Upsilon\left(\frac{k}{2n}\right)\Big[1- \frac{c_\alpha}{4 \min\{ k, 2n-k\} }\Big]~,
\end{align*}
because one of the factors $k$ and $2n-k$ equals $\min\{ k, 2n-k\}$ and
the other one does not exceed $2n-1$. Thus
\begin{align*}
\sum_j \sigma_{k,j} \varphi_n\left(\frac{j}{n}\right) &=
\frac{\theta_\alpha}{n^\alpha} + \frac{1}{n^{\alpha/2}} \sum_j
\sigma_{k,j} \Upsilon\left(\frac{j}{n}\right)\\
&\geq \frac{\theta_\alpha}{n^\alpha} + \frac{1}{n^{\alpha/2}}
\Upsilon\left(\frac{k}{2n}\right) \Big[1 - \frac{c_\alpha}{4 \min(k, 2n-k)}\Big]~.
\end{align*}
We have to compare that with
$$(1+\gamma) \varphi_{2n} \left(\frac{k}{2n}\right) = \Big[\frac{\theta_\alpha}{(2n)^\alpha} +\frac{1}{(2n)^{\alpha/2}}\Upsilon\left(\frac{k}{2n}\right)\Big](1+\gamma)~.$$
Clearly, $\frac{\theta_\alpha}{n^\alpha} -
\frac{(1+\gamma)\theta_\alpha}{(2n)^\alpha} \geq (1-
\frac{1+\gamma}{2^\alpha}) \frac{\theta_\alpha}{n^\alpha} \geq 0$ if
$\gamma < 2^\alpha -1$. Subtracting the second terms, we get
$$\Upsilon\left(\frac{k}{2n}\right)\frac{1}{n^{\alpha/2}} \Big[1- \frac{1+\gamma}{2^{\alpha/2}} - \frac{c_\alpha}{4 \min(k, 2n-k)}\Big]~,$$
which is non-negative if $\gamma < 2^{\alpha/2} -1$ and if $k$ or
$2n-k$ is larger than some constant $K_* =K_*(\alpha, \gamma)$. But,
for $\min(k, 2n-k) \leq K_*$, we have $\Upsilon(\frac{k}{2n}) \leq
[\frac{K_*}{2n}]^{\alpha/2}$ and, thereby, the difference is (in absolute value) at most
$\frac{c_\alpha}{4} \frac{K_*^{\alpha/2}}{n^\alpha}$, which is
dominated by $(1-\frac{1+\gamma}{2^\alpha})
\frac{\theta_\alpha}{n^\alpha}$, provided that $\theta_\alpha$ was chosen
large enough. This proves (\ref{want}).
\medskip

An immediate corollary is that
$$(B_{n/b}\varphi_{n/b}) \succeq_n (1+\gamma)B_n\varphi_n$$
for every $n\in \Lambda \backslash \{n_0\}$. Thus, the Bernstein
coefficients of
$$[x+(1-x)]^r (B_{n/b}\varphi_{n/b}) [x+(1-x)]^{n-n/b} - [x+(1-x)]^r B_n \varphi_n$$
are at least as large as those of $\gamma [x+ (1-x)]^r B_n
\varphi_n$. Since $\varphi_n \geq (\Delta_n)^\alpha$ (provided that
$\theta_\alpha\geq 1$, of course), we see that the latter dominate the
Bernstein coefficients of $Q_{n,r}(f-f_{n/b})$ with some small constant.
Choosing $D$ large enough, we turn this into true domination, which
finishes the proof of ``monotonicity'' of the Bernstein coefficients
of $g_n$ and $h_n$.

\section{Revisiting the claim of Lorentz \label{sec_fedya}}

The goal of this section is to demonstrate that Lorentz' Claim~\ref{res_Lorentz}
made in \cite{Lorentz_paper} is invalid. Our counterexample will be constructed
in several steps. We begin with some elementary observations about Bernstein
polynomials.

\begin{lemma}
Let $\B_n[a,b] \eqbd \{\sum_{k=0}^n c_k (x-a)^k (b-x)^{n-k} : c_k \geq
0 \}$. Then
\begin{description}

\item{(a)} $\B_n[a,b] \subset \B_{n+1} [a, b]$,

\item{(b)} $\B_n[a,b] \cdot \B_m[a, b]\subset \B_{n+m}[a,b]$,

\item{(c)} $\B_n[a,b] \subset \B_n[c,d]$ for every subinterval $[c,d]$ of
the interval $[a,b]$,

\item{(d)} $\B_n[a,b]$ is a convex cone of functions.

\end{description}
\end{lemma}

\proof

(a) Multiply by $1=\frac{1}{b-a} [(x-a) + (b-x)]$ and distribute.

(b) Multiply out.

(c) $x-a = (c-a) + (x-c) \in \B_0[c,d] + \B_1 [c,d] = \B_1 [c,d]$ and
$$b-x = (b-d) + (d-x) \in \B_0[c,d]+ \B_1[c,d] = \B_1[c,d]~. \qquad
\qquad \qquad \qquad \qquad \quad \; $$

Hence,
$(x-a)^k (b-x)^{n-k} \in \B_1[c,d]^n \subset \B_n[c,d]~.$

(d) Obvious.
\eop

\begin{lemma}
Suppose that $p$ is a polynomial of degree $n$ with real
coefficients such that $p(0)>0$ and $p$ has no roots in the unit
disc $\{|z| \leq 1\}$. Then $p\in \B_n [-1, 1]$.
\end{lemma}

\proof
We have $p(x) = \alpha \prod_\beta (x-\beta) \prod_\gamma(\gamma
-x) \prod_\lambda  (x - \lambda) (x-\overline{\lambda})$ where $\beta$ are negative
roots, $\gamma$ are positive roots, $\lambda$ are complex roots
with positive imaginary parts, and $\alpha>0$. Now
$$x-\beta = (x+1) + (-\beta - 1)\quad  \hbox{and} \quad  -\beta -1>0~.$$
Thus, $x-\beta \in \B_1[-1,1]$ for all $\beta$. Similarly, $\gamma
- x \in \B_1[-1,1]$ for all $\gamma$. Now,
$$ 
(x - \lambda)(x-\overline{\lambda}) = x^2 - 2 Re(\lambda x) + |\lambda|^2 \hbox{ is a convex combination of } (|\lambda| - x)^2 \hbox{ \rm and }  (x + |\lambda|)^2 \, .
$$
Moreover, since
$|\lambda| - x  \in \B_1[-1,1]  \hbox{ \rm and } |\lambda| + x  \in \B_1[-1,1] \,, $
we infer that
$ (|\lambda| - x)^2 \in \B_2[-1,1] \,$  and $(|\lambda| + x)^2 \in \B_2[-1,1] \, .$
\eop

\begin{lemma}
The Taylor polynomial $P_{2n}$ of degree $2n$ of the function
$\ee^{-x^2}$ at 0 has no roots in the disc $\{|z| \leq
\frac{\sqrt{n}}{\ee}\}$.
\end{lemma}

\proof
Let $|z| \leq \frac{\sqrt{n}}{\ee}$. Then
\begin{align*}
|\ee^{-z^2} - P_{2n}(z)| \leq \sum_{k>n}
\frac{|z|^{2k}}{k!} \leq
\sum_{k>n}\left(\frac{\ee|z|^2}{k}\right)^k
\leq \sum_{k>n}\left(\frac{\ee|z|^2}{n}\right)^k \leq
\sum_{k>n} \ee^{-k} < \ee^{-n} \leq
|\ee^{-z^2}|~,
\end{align*}
and the result follows.  \eop

In the sequel, we will make use  of the inequality
$$|\ee^{-z^2} - P_{2n}(z)| < \ee^{-n} \quad \hbox{  for  } \; |z|\leq \frac{\sqrt{n}}{\ee}$$
obtained in the course of the last proof.

The following lemma is proved analytically, but the motivation of the construction is probabilistic.
The Bernstein approximation $B_nf$ of a function $f$ can be viewed as the expectation of $f(X/n)$ 
where $X$ is a Binomial random variable with parameters $n$ and $x$. The Central limit theorem yields 
convergence of scaled Binomial variables to Gaussian variables, so the Bernstein approximation is
close to the convolution of $f$ with a suitable Gaussian variable.

\begin{lemma} \label{lem-convolve}
Suppose that $\nu$ is a positive measure on $\R$ such that $g \eqbd 
\nu * \ee^{-nx^2}$ is bounded on the entire real line. Then
there exists $p_n \in \B_{200n}[-1,1]$ such that
$\|g-p_n\|_{L^\infty[-1,1]} \leq 3 \ee^{-n} \|g\|_\infty$.
\end{lemma}

\proof
Note that $x \mapsto P_{200n}(\sqrt{n}(x-t))$ has no roots in the disc
$\{|z-t|\leq \frac{10}{\ee}\}$. If $|t| \leq 2$, this disc
contains the disc $\{|z| \leq 1\}$, so $P_{200n} (\sqrt{n}(x-t))
\in \B_{200n}[-1,1]$. Now put
$$p_n = \nu|_{[-2,2]} * P_{200n}(\sqrt{n}~\cdot)\in \B_{200n} [-1,1]~.$$
For all $x\in [-1,1]$, we have
\begin{align*}
&|g(x) - p_n(x)|= \int_2^\infty \ee^{-n(x-t)^2}~\dd\nu(t) +
\int_{-\infty}^{-2} \ee^{-n(x-t)^2}~\dd\nu(t)\\
& + \int_{-2}^2 |\ee^{-n(x-t)^2} -
P_{200n}(\sqrt{n}(x-t))|~\dd\nu(t) \bdeq I_1 +I_2 +I_3~.
\end{align*}
But $|\ee^{-n(x-t)^2} - P_{200n}(\sqrt{n} (x-t))| \leq
\ee^{-100n}$ as long as $|x-t| \leq \frac{10}{\ee}$.
So $I_3 \leq \ee^{-100n} \int_{-2}^2 \dd\nu(t)$. On the other hand,
since $$\int_{-1}^1\ee^{-nx^2}~\dd x =
\frac{1}{\sqrt{n}}\int_{-\sqrt{n}}^{\sqrt{n}}\ee^{-x^2}~\dd x \geq
\frac{1}{2\sqrt{n}}~,$$
the definition of $g$ implies that
$$\|g\|_\infty \cdot 6 \geq \int_{-3}^3 g(x) \dd x \geq \int_{-1}^1  \ee^{-nx^2} \dd x 
\int_{-2}^2 \dd\nu(t) =  \frac{1}{2\sqrt{n}} \int_{-2}^2 \dd\nu(t)~.$$
So $I_3\leq \ee^{-100n} 12 \sqrt{n} \| g\|_\infty \leq
\ee^{-95n} \|g\|_\infty$.

Now, since for every $y>0, z>1$, we have $\ee^{-n(y+z)^2}
\leq \ee^{-ny^2} \ee^{-n}$, we obtain
\begin{eqnarray*}
I_1 & = & \int_2^\infty \ee^{-n(t-x)^2} ~\dd\nu(t) =
\int_2^\infty \ee^{-n((t-2)+ (2-x))^2} ~\dd\nu(t)\\
&\leq&\int_2^\infty \ee^{-n} \cdot
\ee^{-n(t-2)^2}~\dd\nu(t)\le  \ee^{-n} g(2) \leq
\ee^{-n} \|g\|_\infty~,
\end{eqnarray*}
and, similarly, $I_2 \leq \ee^{-n} \|g\|_\infty$. Bringing
these three estimates together, we arrive at the conclusion of the lemma.
\eop

\begin{corollary}
Let $E_n \eqbd \{ \nu * \ee^{-nx^2}\}$, where the measure $\nu$ satisfies 
the assumption of Lemma~\ref{lem-convolve}. If $f: \R \to [0,1]$
can be approximated by functions $g_n\in E_n$ on the entire
line with an error $O(n^{-\alpha/2})$, then $f$ can be
approximated by $p_n \in \B_{n} [-1,1]$ with an error
$O(n^{-\alpha/2})$ on $[-1,1]$.
\end{corollary}
\proof
Obvious from Lemma \ref{lem-convolve}.
\eop

Now fix $\alpha \in (0,1)$. Our next task will be to construct
a function $f : \R \to [0,1]$ that is approximable by functions
$g_n \in E_{\pi n}$ with an error $O(n^{-\alpha/2})$ but is not
in the class $C^\alpha[-1/2,1/2]$.  Note that $E_\lambda \subset
E_{\lambda'}$ whenever $\lambda < \lambda'$, so it does not matter
whether we consider only integer values or all real values of $n$
in our statement.

Fix $h\in (0,1)$ and $m\in \N$ and define
$$ f_{h,m}(x)  \eqbd h \sum_{k\in \Z} \ee^{-\pi m(x-\frac{kh}{\sqrt{m}})^2} =
h \sum_{k\in \Z} \ee^{-\pi h^2(k-\frac{x\sqrt{m}}h)^2}. $$
Recalling that the Fourier transform of the function
$x \mapsto h\ee^{-\pi h^2 x^2}$ is $y \mapsto \ee^{-\pi y^2/h^2}$
and using the Poisson summation formula
$$  \sum_{k\in \Z} F(k+x) = \sum_{\ell\in \Z} \widehat{F}(\ell) \ee^{2\pi i \ell x},    $$
we get
$$ f_{h,m} = \sum_{\ell\in \Z} \ee^{-\pi \ell^2 / h^2}
\ee^{-2\pi i \ell x \sqrt{m} /h}. $$
This representation immediately implies that $f_{h,m}$ attains its
maximum at $x=0$, and that
\begin{eqnarray*}
 | f_{h,m} -1 | & \leq & \sum_{\ell \in \Z \setminus \{0\}}
\ee^{-\pi \ell^2/h^2} \leq 2 \ee^{-\pi/h^2}  (1+ \sum_{\ell \geq 3} \ee^{-\pi \ell / h^2}  )
\leq 4 \ee^{-\pi / h^2} \\
f_{h,m}(0) - f_{h,m} \left( {h \over 2\sqrt{m} }  \right) &
= & 2 \sum_{\ell \; {\rm odd}} \ee^{-\pi \ell^2 / h^2}
= 2 \sum_{|\ell|=1}  \ee^{-\pi \ell^2 / h^2} \; + \cdots
  \geq 4 \ee^{-\pi/h^2}.
\end{eqnarray*}
Also note that $f_{h,m}\in E_{\pi m}.$

Now let $\Lambda$ denote the set $\{ 2^j : j=2,3,4,\ldots \}$. Choose $h_m$
so that   $$ { \ee^{-\pi/h_m^2} \over (\log_2 m)^2 } = {1\over m^{\alpha/2} }.$$
 This choice makes sense for $m \geq m_0(\alpha) \geq 4$.
Define $\Lambda' \eqbd \{ m \in \Lambda : m \geq m_0(\alpha) \} $ and
$$ f \eqbd \sum_{m\in \Lambda'} {1\over (\log_2 m)^2} f_{h_m,m}. $$
For every $n\in \N$, let
$$ g_n \eqbd \sum_{m\in \Lambda', m\leq n } {1\over (\log_2 m)^2} f_{h_m,m}
+ \sum_{m\in \Lambda', m>n } {1\over (\log_2 m)^2} \in E_{\pi n}, $$
since $g_n$ is a finite sum of elements of $E_{\pi m}$ with $m\le n$
plus a constant.  Now,
\begin{eqnarray*}
|f-g_n| & \leq  & \sum_{m\in \Lambda', m > n } {1\over (\log_2 m)^2}
\|f_{h_m,m}-1 \|_\infty \leq
 \sum_{m\in \Lambda', m > n } {1\over (\log_2 m)^2} 4 \ee^{-\pi/h_m^2}  \\
& = &
4 \sum_{m\in \Lambda', m > n }  m^{-\alpha/2} \leq \const n^{-\alpha/2}.
\end{eqnarray*}
On the other hand, for every $m\in \Lambda'$, we have
$$ f(0) - f\left( {h_m \over 2\sqrt{m}}   \right)  \geq \left(2 \sum_{\ell \;{\rm odd}}
\ee^{-\pi \ell^2/h_m} \right)  (\log_2 m)^{-2} > 4 { \ee^{-\pi/h_m^2}\over (\log_2 m)^2}= 4 m^{-\alpha/2}.  $$
Thus, $$ \| f \|_{C^\alpha[-1/2,1/2]} \geq {4 m^{-\alpha/2} \over
\Bigl(h_m /(2\sqrt{m})\Bigr)^\alpha} = 16 h_m^{-\alpha} \to \infty   \quad {\rm as}
\; m\to \infty, $$
so $f$ is not in the class $C^\alpha[-1/2,1/2]$. So, we have obtained
a function $f\notin C^\alpha[-1/2,1/2]$ such that it can be approximated
by polynomials $p_n \in \B_n [-1,1]$ at the rate $O(n^{-\alpha/2})$.

Consider the function $\widetilde{f}(x) \eqbd f(x) \cdot x (1-x)$
and the polynomials $\widetilde{p}_n (x) \eqbd p_n(x) \cdot x (1-x)$.
The polynomials are in $\B_{n+2}[-1,1]$ and the function
$\widetilde{f}$ satisfies the condition
$$  | \widetilde{f}(x) - \widetilde{p}_n (x)  | \leq x (1-x)
\const n^{-\alpha/2} \leq \const \left( \sqrt{x(1-x)\over n} \right)^\alpha 
\le \const  (\Delta_n(x))^\alpha. $$ Claim \ref{res_Lorentz} is thus disproved.

\section{Further questions and remarks}  \label{sec_moments}
In this section we will make a few additional remarks on this and some
related problems.
We begin by discussing a conjectural characterization of simulation rates
in case $\alpha$ is an integer. In that case, the classical problem of
approximating a given function by polynomials of degree at most $n$
already has a somewhat different solution, as we now explain.

\begin{definition} Let $\alpha\in \N$.
A function $f$ is said to be in the {\bf\em smoothness class\/}
${C^\alpha}^*[0,1]$  if $f$ is  $r\eqbd  \alpha{-}1$  times differentiable
and the following condition holds:

The {\bf\em symmetric modulus of continuity\/} of
$f^{(r)}$
$$\omega^*(f^{(r)},h) \eqbd \sup_{t<h,\; x\in [t,1-t]
 } | f^{(r)}(x+t)-2 f^{(r)}(x) + f^{(r)}(x-t)| $$
is of order $O(h)$. In that case, we use the
notation $$ \|f \|_{{C^\alpha}^*} \eqbd \sup_h  {\omega^*(f^{(r)},h) \over h  }.$$
\end{definition}

\begin{remark} The class ${C^\alpha}^*[0,1]$ is also known as the
{\bf\em generalized Lipschitz class,\/} and, for $\alpha=1$,
as the {\bf\em Zygmund class\/} we already defined in Section~\ref{sec-necessity}
 (also see \cite[Chap.2, Sec.9]{DeVoreL}).
\end{remark}

The characterization of the best polynomial approximation in case $\alpha\in \N$
is then given by the following result.

\begin{result}[{see, e.g.,~\cite[Chapter 8,~Theorem~6.3]{DeVoreL}}] \label{res_best*}
Let $\alpha \in \N$. There exists a sequence of polynomials $(p_n)$, where the
 degree of $p_n$ is at most $n$, satisfying
$$ | p_n(x) -  f(x)| = O((\Delta_n(x))^{2\alpha}) \qquad \hbox{\rm for all}
\;\; x\in [0,1] $$ if and only if $f\in {C^\alpha}^*[0,1]$.
\end{result}

Motivated by this result on unrestricted polynomial approximation, we therefore conjecture
a corresponding characterization of simulation rates.

\begin{conjecture}  Let $\alpha \in \N$.
Let $f\in {C^\alpha}^* [0,1]$ be a function bounded strictly between $0$ and $1$.
Then $f$ can be simulated at the rate $(\Delta_n(x))^\alpha$ on $[0,1]$. Precisely,
there exist polynomials $g_n$ and $f_n$ satisfying conditions (i) -- (iv) of
Result~\ref{res_reduction} and bound~(\ref{app_order}).
\end{conjecture}

In Theorem~\ref{thm-int} of Section~\ref{sec-necessity}, we have already
verified the converse: if $f$ is simulable at the rate $(\Delta_n(x))^\alpha$ on the interval
$[0,1]$ where $\alpha \in \N$,  then $f\in {C^\alpha}^* [0,1]$.

Finally, we note that that for any $\alpha>0$, it is natural to ask   which functions $f$    
can be simulated with a finite $\alpha$ moment, i.e., when is there a simulation algorithm 
for an $f(p)$-coin such that the number $N$ of tosses of $p$-coins  and fair coins it uses 
satisfies $${\mathbf E}_p(N^\alpha) \; ( = \sum_{n=1}^\infty  n^\alpha \P_p (N>n) )\; <\; \infty.$$
We suspect that the precise criterion should involve the Besov smoothness of $f$, with proper 
attention to boundary effects; see, e.g.,~\cite[pp.\ 54--57]{DeVoreL} for the definition and basic 
properties of Besov spaces.

\newpage

\bibliographystyle{plain}
\bibliography{bernstein}
\end{document}